\documentclass[12pt, reqno]{amsart}
\usepackage{amssymb,verbatim,amscd,amsmath,graphicx}
\usepackage{graphicx}
\usepackage{caption}
\usepackage{subcaption}
\usepackage{pdfsync}
\usepackage{fullpage}
\usepackage{color}
\usepackage{enumitem}
\usepackage{tikz}
\usetikzlibrary{patterns}
\usetikzlibrary{calc}
\usetikzlibrary{matrix}
\usepgflibrary{arrows}
\usetikzlibrary{arrows}
\usetikzlibrary{decorations.pathreplacing}
\usetikzlibrary{decorations.pathmorphing}
\numberwithin{equation}{section}
\newtheorem{theorem}{Theorem}[section]
\newtheorem{lemma}[theorem]{Lemma}

\newtheorem{corollary}[theorem]{Corollary}
\newtheorem{conjecture}[theorem]{Conjecture}

\theoremstyle{definition}
\newtheorem{definition}[theorem]{Definition}
\newtheorem{def-prop}[theorem]{Definition-Proposition}
\newtheorem{remark}[theorem]{Remark}
\newtheorem{example}[theorem]{Example}

\newtheorem*{acknowledgement}{Acknowledgement}

\newtheorem{question}[theorem]{Question}
\newtheorem{problem}[theorem]{Problem}
\newtheorem*{Mysketch}{Sketch of proof} % \newtheorem establishes the object heading
\newenvironment{sketch}    % this is the environment name for the input
  {\pushQED{\qed}\begin{Mysketch}}
  {\popQED\end{Mysketch}}

\DeclareMathOperator{\Tor}{Tor}
\DeclareMathOperator{\reg}{reg}

\DeclareMathOperator{\adeg}{adeg}

\DeclareMathOperator{\link}{link}
\DeclareMathOperator{\del}{del}
\DeclareMathOperator{\depth}{depth}
\DeclareMathOperator{\pd}{pd}
\DeclareMathOperator{\supp}{supp}

\DeclareMathOperator{\pol}{pol}

\newcommand{\ZZ}{{\mathbb Z}}
\newcommand{\NN}{{\mathbb N}}

\def\mm{{\frak m}}

\def\H{{\mathcal H}}

\def\1{{\bf 1}}
\def\0{{\bf 0}}

\begin{document}

\title{Regularity of edge ideals and their powers}

\author{Arindam Banerjee}
\address{Ramakrishna Mission Vivekenanda University, Belur, West Bengal, India}
\email{123.arindam@gmail.com}
\urladdr{https://http://maths.rkmvu.ac.in/~arindamb/}

\author{Selvi Beyarslan}
\address{University of South Alabama, Department of Mathematics and Statistics, 411 University Boulevard North, Mobile, AL 36688-0002, USA}
\email{selvi@southalabama.edu}
\urladdr{https://selvikara.wordpress.com/selvi/}

\author{Huy T\`ai H\`a}
\address{Tulane University \\ Department of Mathematics \\
6823 St. Charles Ave. \\ New Orleans, LA 70118, USA}
\email{tha@tulane.edu}
\urladdr{http://www.math.tulane.edu/$\sim$tai/}

\begin{abstract} We survey recent studies on the Castelnuovo-Mumford regularity of edge ideals of graphs and their powers. Our focus is on bounds and exact values of $\reg I(G)$ and the asymptotic linear function $\reg I(G)^q$, for $q \ge 1$, in terms of combinatorial data of the given graph $G$.
\end{abstract}

\maketitle

%%%%%%%%%%%%%%%%%%%%%%%%%%%%%%%%%%%%%%%%%%%%%%

\section{Introduction}

Monomial ideals are classical objects that live at the crossroad of three areas in mathematics: algebra, combinatorics and topology. Investigating monomial ideals has led to many important results in these areas. The new construction of edge ideals of (hyper)graphs has again resurrected much interest in and regenerated a large amount of work on this class of ideals (cf. \cite{HVTsurvey, MV, V} and references therein). In this paper, we survey recent works on the Castelnuovo-Mumford regularity of edge ideals of graphs and their powers.

Castelnuovo-Mumford regularity is an important algebraic invariant which, roughly speaking, measures the complexity of ideals and modules. Restricting to the class of edge ideals of graphs, our focus is on studies that relate this algebraic invariant to the combinatorial data of given graphs. Our interest in powers of edge ideals is driven by the well-celebrated result of Cutkosky, Herzog and Trung \cite{CHT}, and independently, Kodiyalam \cite{Ko}, that for any homogeneous ideal $I$ in a standard graded $k$-algebra $R$, the regularity of $I^q$ is asymptotically a linear function in $q$; that is, there exist constants $a$ and $b$ such that for all $q \gg 0$, $\reg I^q = aq+b$. Generally, the problem of finding the exact linear form $aq+b$ and the smallest value $q_0$ such that $\reg I^q = aq+b$ for all $q \ge q_0$ has proved to be very difficult. We, thus, focus our attention also on the problem of understanding the linear form $aq+b$ and the value $q_0$ for edge ideals via combinatorial data of given graphs.

The on-going research program in which algebraic invariants and properties of edge ideals and their powers are investigated through combinatorial structures of corresponding graphs has produced many exciting results and, at the same time, opened many further interesting questions and conjectures. It is our hope to collect these works together in a systematic way to give a better overall picture of the problems and the current state of the art of this research area. For this purpose, we shall state theorems and present sketches of the proofs; instead of giving full detailed arguments, our aim is to exhibit general ideas behind these results, the similarities and differences when developing from one theorem to the next.

This paper can be viewed as a complement to the survey done in \cite{HaSurvey}. While in \cite{HaSurvey} the focus was on the regularity of squarefree monomial ideals in general, our attention in this paper is restricted mostly to edge ideals of graphs, but we enlarge our scope by discussing also the regularity of powers of edge ideals. The later is an important area of study itself, with a deep motivation from geometry, and has seen a surge of interest in the last few years.

The paper is outlined as follows. In Section \ref{chapter2}, we collect notations and terminology used in the paper. In Section \ref{chapter3}, we present necessary tools which were used in works in this area. Particularly, we shall give Hochster's and Takayama's formulas, which relate the graded Betti numbers of a monomial ideal to the reduced homology groups of certain simplicial complexes. We shall also describe inductive techniques that have been the backbones of most of the studies being surveyed. In Section \ref{chapter4}, we survey results on the regularity of edge ideals of graphs. This section is divided into two subsections; the first one focusses on bounds on the regularity in terms of combinatorial data of the graph, and the second one exhausts cases where the regularity of an edge ideal can be computed explicitly. Section \ref{chapter5} discusses the regularity of powers of edge ideals. This section again is divided into two subsections, which in turn examine bounds for the asymptotic linear function of the regularity of powers of edge ideals and cases when this asymptotic linear function can be explicitly described. In Section \ref{chapter6}, we recall a number of results extending the study of edge ideals of graphs to hypergraphs. Since, our focus in this paper is on mostly edge ideals of graphs, our results in Section \ref{chapter6} will be representative rather than exhaustive. We end the paper with Section \ref{chapter7}, in which we state a number of open problems and questions which we would like to see answered. We hope that these problems and questions would stimulate further studies in this research area.

\begin{acknowledgement} The authors would like to thank the organizers of SARC (Southern Regional Algebra Conference) 2017 for their encouragement, which led us to writing this survey. The last named author is partially supported by the Simons Foundation (grant \#279786) and Louisiana Board of Regent (grant \#LEQSF(2017-19)-ENH-TR-25).
\end{acknowledgement}

\section{Preliminaries}\label{chapter2}

   In this section we recall preliminary notations and terminology of combinatorics and algebra that will be used throughout this survey. In Subsection \ref{sub.comb} we give various definitions on graphs, hypergraphs, and simplicial complexes. In Subsection \ref{sub.alg} we recall basic homological algebra terminology. Finally, in the last subsection we define various algebraic objects associated to (hyper)graphs and simplicial complexes.
   
\subsection{Combinatorial preliminaries} \label{sub.comb}
For any finite simple graph $G$, with set of vertices $V(G)$ and set of edges $E(G)$, we define some graph-theoretic notions as follows.

\begin{figure}[h]
\begin{tikzpicture}
[scale=.8, vertices/.style={draw, fill=black, circle, inner sep=1.3pt}]
\useasboundingbox (-2.5,-2.4) rectangle (2.5,2.4);
\node [anchor=base] at (-3.3,-.1){$H'$:};
\node [vertices] (1) at (0:2) {};
\node [anchor=base] at (2.4,-.1) {$w$};
\node [vertices] (2) at (60:2) {};
\node [anchor=base] at (60:2.3) {$z$};
\node [vertices] (3) at (120:2) {};
\node [anchor=base] at (120:2.3) {$y$};
\node [vertices] (4) at (180:2) {};
\node [anchor=base] at (-2.4, -.1) {$x$};
%\node [vertices] (5) at (240:2) {};
%\node [anchor=base] at (240:2.5) {$t$};
\node [vertices] (6) at (300:2) {};
\node [anchor=base] at (300:2.5) {$s$};
\foreach \to/\from in {2/1, 3/1, 4/1, 4/3,6/1}
	\draw [-] (\to)--(\from);
\end{tikzpicture}
\hspace*{2.5cm}\begin{tikzpicture}
[scale=.8, vertices/.style={draw, fill=black, circle, inner sep=1.3pt}]
\useasboundingbox (-2.5,-2.4) rectangle (2.5,2.4);
\node [anchor=base] at (-3.3,-.1){$H''$:};
\node [vertices] (1) at (0:2) {};
\node [anchor=base] at (2.4,-.1) {$w$};
\node [vertices] (2) at (60:2) {};
\node [anchor=base] at (60:2.3) {$z$};
\node [vertices] (3) at (120:2) {};
\node [anchor=base] at (120:2.3) {$y$};
\node [vertices] (4) at (180:2) {};
\node [anchor=base] at (-2.4, -.1) {$x$};
%\node [vertices] (5) at (240:2) {};
%\node [anchor=base] at (240:2.5) {$t$};
\node [vertices] (6) at (300:2) {};
\node [anchor=base] at (300:2.5) {$s$};
\foreach \to/\from in {2/1, 4/1, 4/3,6/1}
	\draw [-] (\to)--(\from);
\end{tikzpicture}
\caption{Induced subgraphs} \label{fig:inducedsubgraph}
\end{figure}
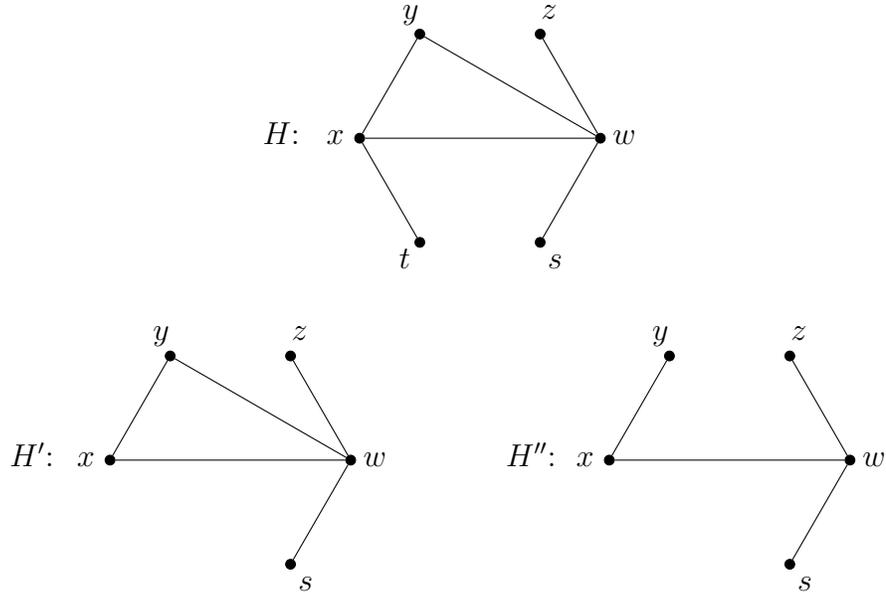

For any $A\subseteq V(G)$, the\textit{ induced subgraph} on $A$ is the graph whose set of vertices is $A$ and whose edges are exactly the edges of $G$ that join two elements of $A$. In the following example (Figure \ref{fig:inducedsubgraph}) $H'$ is an induced subgraph of $H$, while $H''$ is not an induced subgraph of $H$ as it misses the edge $yw$ of $H$.

\begin{figure}[h]
\centering
\begin{tikzpicture}
[scale=.8, vertices/.style={draw, fill=black, circle, inner sep=1.3pt}]
\useasboundingbox (-2.5,-2.4) rectangle (2.5,2.4);
\node [anchor=base] at (-3.3,-.1){$H$:};
\node [vertices] (1) at (0:2) {};
\node [anchor=base] at (2.4,-.1) {$w$};
\node [vertices] (2) at (60:2) {};
\node [anchor=base] at (60:2.3) {$z$};
\node [vertices] (3) at (120:2) {};
\node [anchor=base] at (120:2.3) {$y$};
\node [vertices] (4) at (180:2) {};
\node [anchor=base] at (-2.4, -.1) {$x$};
\node [vertices] (5) at (240:2) {};
\node [anchor=base] at (240:2.5) {$t$};
\node [vertices] (6) at (300:2) {};
\node [anchor=base] at (300:2.5) {$s$};
\foreach \to/\from in {2/1, 3/1, 4/1, 4/3, 5/4, 6/1}
	\draw [-] (\to)--(\from);
\end{tikzpicture}
%\caption{Example of a graph} \label{fig:agraph}
\end{figure}

For any vertex $x$ in a graph G the \textit{degree} of $x$, denoted by $d(x)$ is the number of vertices connected to $x$. For any vertex $x$ in a graph $G$  the \textit{neighborhood of} $x$, denoted by $N_G[x]$ is the set consisting of $x$ and all its neighbors. For the graph $H''$ above,  $d(x)=2$ and $N_G[x]=\{x,y,w\}$. For any graph $G$ and any vertex $x$, by $G\setminus x$ we denote the induced subgraph on $V(G)\setminus \{x\}$.

For any graph $G$, the \textit{complement graph}, denoted by $G^c$, is the graph whose vertices are the vertices of $G$ and whose edges are the non-edges of $G$, i.e., for $a,b \in V(G), ab$ is an edge in $G^c$ if and only if $ab$ is \textit{not} an edge in $G$. For the graph $H$ in Figure \ref{fig:inducedsubgraph}, the following is $H^c$.

\begin{figure}[h]
\begin{tikzpicture}
[scale=.8, vertices/.style={draw, fill=black, circle, inner sep=1.3pt}]
\useasboundingbox (-2.5,-2.4) rectangle (2.5,2.4);
\node [anchor=base] at (-3.3,-.1){$H^c$:};
\node [vertices] (1) at (0:2) {};
\node [anchor=base] at (2.4,-.1) {$w$};
\node [vertices] (2) at (60:2) {};
\node [anchor=base] at (60:2.3) {$z$};
\node [vertices] (3) at (120:2) {};
\node [anchor=base] at (120:2.3) {$y$};
\node [vertices] (4) at (180:2) {};
\node [anchor=base] at (-2.4, -.1) {$x$};
\node [vertices] (5) at (240:2) {};
\node [anchor=base] at (240:2.5) {$t$};
\node [vertices] (6) at (300:2) {};
\node [anchor=base] at (300:2.5) {$s$};
\foreach \to/\from in {1/5,2/3,2/4,2/5,2/6,3/5,3/6,4/6,5/6}
	\draw [-] (\to)--(\from);
\end{tikzpicture}
\caption{Complement graph}\label{fig:complementgraph}
\end{figure}
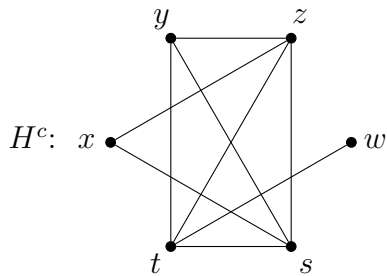

A \textit{cycle} of length $n$  in a graph G is a closed walk along its edges, $x_1x_2,x_2x_3, \ldots, x_{n-1} x_n, x_nx_1$, such that $x_i \neq x_j$ for $ i\neq j$.  We denote the cycle on $n$ vertices by $C_n.$ A \textit{chord} in the cycle $C_n$ is an edge $x_i x_j$ where $x_j \neq x_{i-1},x_{i+1}$. A graph is said to be \textit{chordal} if for any cycle of length greater than  or equal to $4$ there is a chord. A graph is said to be \textit{co-chordal} if the complement of $G$ is chordal. The graph $H$ in Figure \ref{fig:complementgraph} is chordal and $H^c$ is co-chordal.

A \textit{forest} is a graph without any cycles. A \textit{tree} is a connected forest.

A \textit{complete} graph (or \emph{clique}) on $n$ vertices is a graph where for any two vertices there is an edge joining them. It is denoted by $K_n$. 

A \textit{bipartite graph} is a graph whose vertices can be split into two groups such that there is no edge between vertices of same group; only edges are between vertices coming from different groups. It is easy to see that a graph is bipartite if and only if it is without any cycle of odd length. We use $K_{m, n}$ to denote the complete bipartite graph with $m$ vertices on one side, and $n$ on the other. We observe that $H''$ in Figure \ref{fig:inducedsubgraph} has no cycle, so it is a tree.  As a consequence it is also a bipartite graph with bipartition $\{x,z,s\}$ and $\{y,w\}$.

Let $G$ be a graph.  We say two disjoint edges $uv$ and $xy$ form a \emph{gap} in $G$ if $G$ does not have an edge with one
endpoint in $\{u,v\}$ and the other in $\{x,y\}$. A graph without an induced gap is called \emph{gap-free}.  Equivalently, $G$ is gap-free if and only if $G^c$ contains no induced $C_4$. The graph $H$ in Figure \ref{fig:inducedsubgraph} is gap-free. In the following graph $\{x,y\}$ and $\{z,w\}$ form a gap.

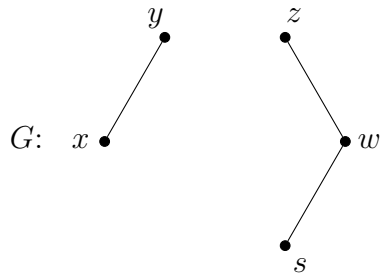
\begin{figure}[h]
\begin{tikzpicture}
[scale=.8, vertices/.style={draw, fill=black, circle, inner sep=1.3pt}]
\useasboundingbox (-2.5,-2.4) rectangle (2.5,2.4);
\node [anchor=base] at (-3.3,-.1){$G$:};
\node [vertices] (1) at (0:2) {};
\node [anchor=base] at (2.4,-.1) {$w$};
\node [vertices] (2) at (60:2) {};
\node [anchor=base] at (60:2.3) {$z$};
\node [vertices] (3) at (120:2) {};
\node [anchor=base] at (120:2.3) {$y$};
\node [vertices] (4) at (180:2) {};
\node [anchor=base] at (-2.4, -.1) {$x$};
%\node [vertices] (5) at (240:2) {};
%\node [anchor=base] at (240:2.5) {$t$};
\node [vertices] (6) at (300:2) {};
\node [anchor=base] at (300:2.5) {$s$};
\foreach \to/\from in {2/1,4/3,6/1}
	\draw [-] (\to)--(\from);
\end{tikzpicture}
\caption{A graph which is not gap-free} \label{fig:notgapfree}
\end{figure}

A \textit{matching} in a graph is a set of pairwise disjoint edges. A matching is called an \textit{induced matching} if the induced subgraph on the vertices of the edges forming the matching has no other edge. We observe that  $\{x,y\}$ and $\{z,w\}$ forms an induced matching for graph $G$ in Figure \ref{fig:notgapfree}.

Any graph isomorphic to $K_{1, 3}$ is called a \emph{claw}. Any graph isomorphic to $K_{1,n}$ is called an \emph{n-claw}. If $n >1$, the vertex with degree $n$ is called the root in $K_{1,n}$.
A graph without an induced claw is called \emph{claw-free}. A graph without an induced \emph{n-claw} is called \emph{n-claw-free}.

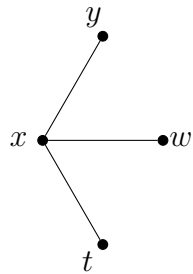
\begin{figure}[h]
\begin{tikzpicture}
[scale=.8, vertices/.style={draw, fill=black, circle, inner sep=1.3pt}]
\useasboundingbox (-2.5,-2.4) rectangle (2.5,2.4);
%\node [anchor=base] at (-3.3,-.1){$G$:};
\node [vertices] (1) at (0:0) {};
\node [anchor=base] at (0.3,-.1) {$w$};
%\node [vertices] (2) at (60:2) {};
%\node [anchor=base] at (60:2.3) {$z$};
\node [vertices] (3) at (120:2) {};
\node [anchor=base] at (120:2.3) {$y$};
\node [vertices] (4) at (180:2) {};
\node [anchor=base] at (-2.4, -.1) {$x$};
\node [vertices] (5) at (240:2) {};
\node [anchor=base] at (240:2.5) {$t$};
%\node [vertices] (6) at (300:2) {};
%\node [anchor=base] at (300:2.5) {$s$};
\foreach \to/\from in {4/5,4/3,4/1}
	\draw [-] (\to)--(\from);
\end{tikzpicture}
\caption{A 3-claw (or simply, a claw)} \label{fig:claw}
\end{figure}

A graph is called a \emph{diamond} if it is isomorphic to the graph with vertex set $\{a,b,c,d\}$ and edge set $\{ab,bc,ac,ad,cd\}$. A graph without an induced diamond is called \emph{diamond-free}.

\begin{figure}[h]
\begin{tikzpicture}
[scale=.8, vertices/.style={draw, fill=black, circle, inner sep=1.3pt}]
\useasboundingbox (-2.5,-2.4) rectangle (2.5,2.4);
%\node [anchor=base] at (-3.3,-.1){$G$:};
%\node [vertices] (1) at (0:2) {};
%\node [anchor=base] at (2.4,-.1) {$w$};
\node [vertices] (2) at (60:2) {};
\node [anchor=base] at (60:2.3) {$z$};
\node [vertices] (3) at (120:2) {};
\node [anchor=base] at (120:2.3) {$y$};
\node [vertices] (4) at (180:2) {};
\node [anchor=base] at (-2.4, -.1) {$x$};
\node [vertices] (5) at (240:2) {};
\node [anchor=base] at (240:2.5) {$t$};
%\node [vertices] (6) at (300:2) {};
%\node [anchor=base] at (300:2.5) {$s$};
%\node [vertices] (7) at (220:2.2) {};
%\node [anchor=base] at (220:2.6) {$y'$};
\foreach \to/\from in { 3/2, 4/3, 5/2, 5/3, 5/4}
	\draw [-] (\to)--(\from);
\end{tikzpicture}
\caption{A diamond} \label{fig:diamond}
\end{figure}
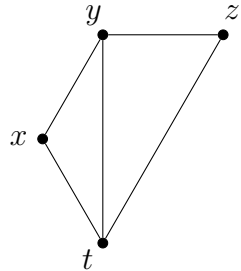

A graph is said to be \textit{planar} if as a $1$-dimensional topological space it can be embedded in the complex plane, i.e., if it can be drawn in the plane in such a way that no pair of edges cross. The following graph in Figure \ref{fig:diamond} is a diamond. We observe that it is also planar.

Any graph isomorphic to the graph with set of vertices $\{w_1,w_2,w_3,w_4,w_5\}$ and set of edges $\{w_1w_3,w_2w_3,w_3w_4,w_3w_5,w_4w_5\}$ is called a \textit{cricket}. A graph without an induced cricket is called \emph{cricket-free}. It is easy to see that a claw-free graph is also cricket-free.

  \begin{figure}[h]
  \includegraphics[width=0.18\linewidth]{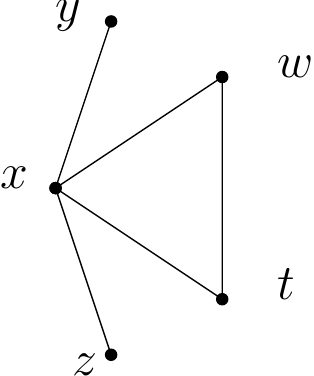}
  \caption{A cricket}
  \label{fig:cricket}
\end{figure}

An edge in a graph is called a \textit{whisker} if one of its vertices has degree one. A graph is called an \textit{anticycle} if its complement is a cycle.

\begin{figure}[h]
  \includegraphics[width=0.7\linewidth]{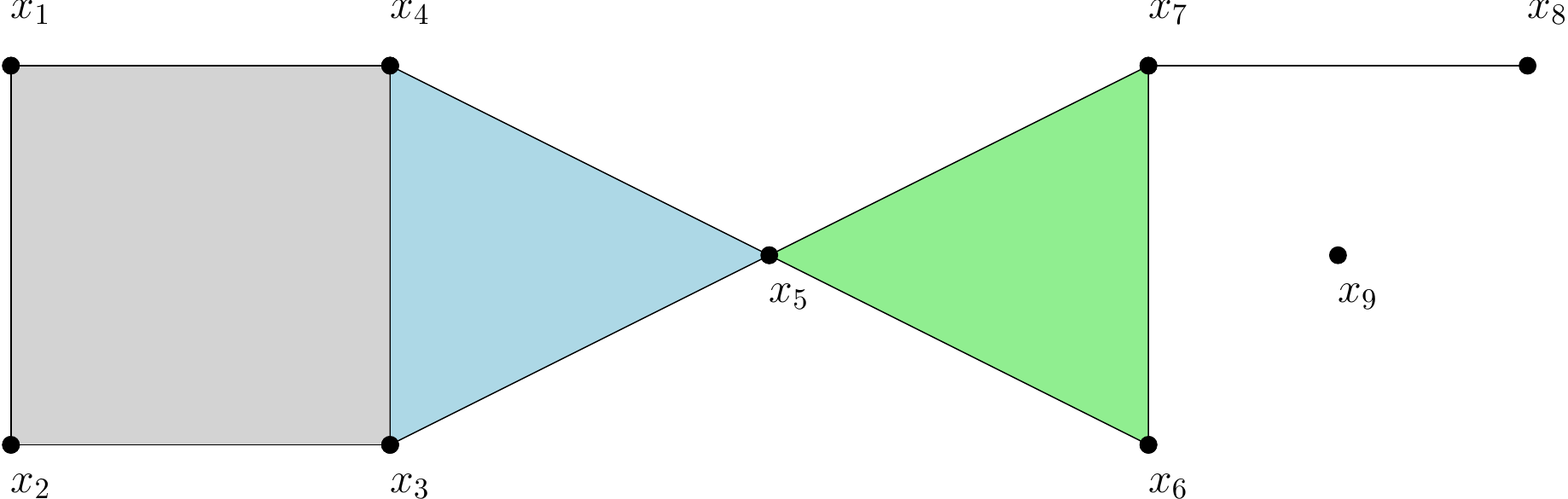}
  \caption{A hypergraph}
  \label{fig:hypergraph}
\end{figure}

A \textit{hypergraph} is the natural higher degree analogue of graphs, in the sense that we allow an edge to be a collection of any number of vertices. Formally speaking, a hypergraph $H = (V(H),E(H))$ consists of the \textit{vertices} $V(H)=\{x_1,\ldots,x_n\}$ and the \textit{edges} $E(H)$, where $E(H)$ is a collection of nonempty subsets of the vertices.  With this notation, a graph $G$ is a hypergraph whose edges are subsets of cardinality 2. By abusing notation, we often identify an edge $\{x_{i_1}, \dots, x_{i_r}\} \in E(H)$ with the monomial $x_{i_1} \dots x_{i_r}$. Figure \ref{fig:hypergraph} depicts a hypergraph with edges $\{x_1x_2x_3x_4, x_3x_4x_5, x_5x_6x_7, x_7x_8,x_9\}.$ In this example, the vertex $x_9$ can also be viewed as an \emph{isolated vertex}.

A \textit{clutter} is a hypergraph none of whose edges contains any other edge as a subset. The hypergraph in Figure \ref{fig:hypergraph} is a clutter.

\begin{figure}[h]
  \includegraphics[width=0.4\linewidth]{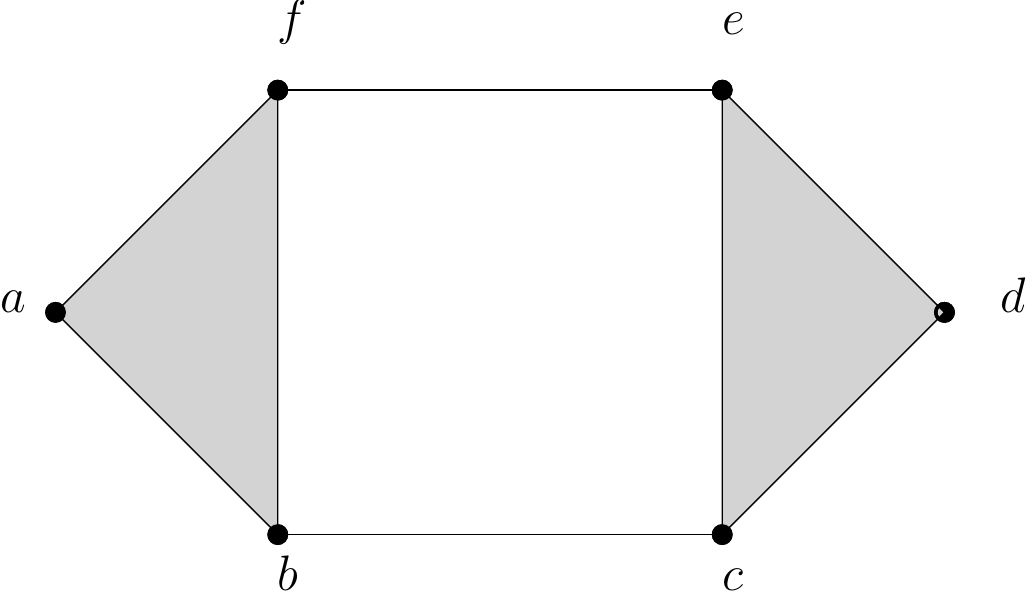}
  \caption{A simplicial complex}
  \label{fig:simplicial}
\end{figure}

A \textit{simplicial complex} with vertices $\{x_1,\ldots,x_n\}$ is subset of the power set of $\{x_1,\ldots,x_n\}$ which is closed under  the subset operation.  The sets that constitute a simplicial complex are called its \textit{faces}. Maximal faces under inclusion are called \textit{facets}. Figure \ref{fig:simplicial} gives a simpicial complex with facets $\{a,b,f\}$, $\{c,d,e\}$, $\{f,e\}$,$\{b,c\}$.

An \textit{independent set} in a graph $G$ is a set of vertices no two of which forms an edge. The independence complex of a graph $G,$ denoted by $\Delta(G),$ is the simplicial complex whose faces are independent sets in $G.$

\begin{figure}[hb]
  \includegraphics[width=0.38\linewidth]{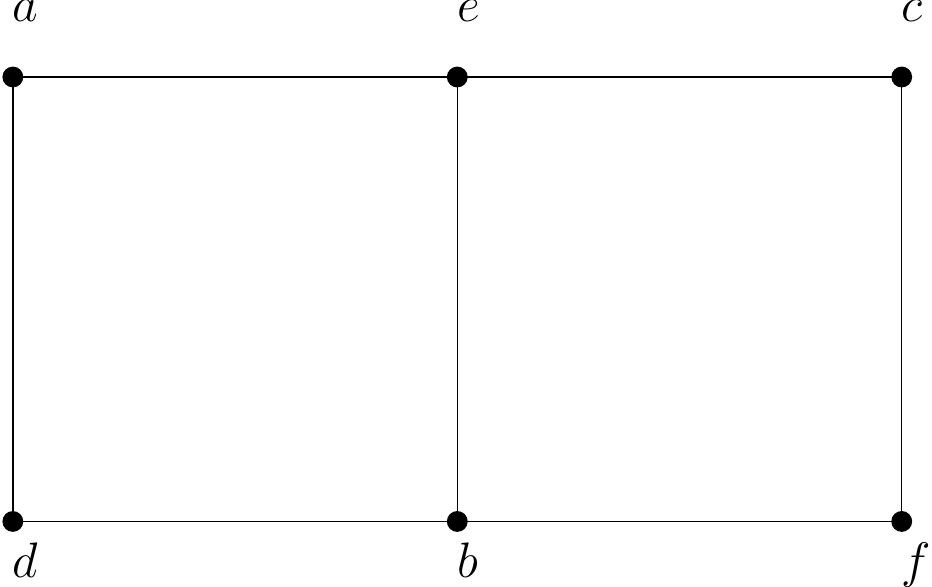}
  \caption{A simple graph whose independence complex is in Figure \ref{fig:simplicial}}
  \label{fig:independence}
\end{figure}

Let $\Delta$ be a simplicial complex, and let $\sigma \in \Delta$.
The \emph{deletion} of $\sigma$ in $\Delta$, denoted by $\del_\Delta(\sigma),$ is the simplicial complex obtained by removing $\sigma$ and all faces containing $\sigma$ from $\Delta$.
The \emph{link} of $\sigma$ in $\Delta$, denoted by $\link_\Delta(\sigma)$, is the simplicial complex whose faces are
$$\{ F \in \Delta ~|~ F \cap \sigma = \emptyset, \sigma \cup F \in \Delta\}.$$

A simplicial complex $\Delta$ is recursively defined to be \emph{vertex decomposable} if either:
\begin{enumerate}
\item[(i)] $\Delta$ is a simplex; or
\item[(ii)] there is a vertex $v$ in $\Delta$ such that both $\link_\Delta(v)$ and $\del_\Delta(v)$ are vertex decomposable, and all facets of $\del_\Delta(v)$ are facets of $\Delta$.
\end{enumerate}
A vertex satisfying condition (2) is called a \emph{shedding vertex}, and the recursive choice of shedding vertices is called a \emph{ shedding order} of $\Delta$.

Recall that a simplicial complex $\Delta$ is said to be \emph{shellable} if there exists a linear order of its facets $F_1, F_2, \dots, F_t$ such that for all $k = 2, \dots, t$, the subcomplex $\Big(\bigcup_{i=1}^{k-1} \overline{F_i}\Big) \bigcap \overline{F_k}$ is pure and of dimension $(\dim F_k - 1)$. Here $\overline{F}$ represents the simplex over the vertices of $F$. It is a celebrated fact that \emph{pure} shellable complexes give rise to \emph{Cohen-Macaulay Stanley-Reisner rings}. Recall also that a ring or module is \emph{sequentially Cohen-Macaulay} if it has a filtration in which the factors are Cohen-Macaulay and their dimensions are increasing. This property corresponds to (\emph{nonpure}) shellability in general.
Vertex decomposability can be thought of as a combinatorial criterion for shellability and sequentially Cohen-Macaulayness. In particular, for a simplicial complex $\Delta$,
$$\Delta \textrm{ vertex decomposable } \Rightarrow \Delta \textrm{ shellable } \Rightarrow \Delta \textrm{ sequentially Cohen-Macaulay}.$$

\subsection{Algebraic preliminaries} \label{sub.alg}
Let $S = K[x_1, \dots, x_n]$ be a polynomial ring over a field $K$.  Let $M$ be a finitely generated $\ZZ^n$-graded $S$-module. It is known that $M$ can be successively approximated by free modules. Formally speaking, there exists an exact sequence of minimal possible length, called \textit{a minimal free resolution} of $M$:
 $$0\longrightarrow \mathbb{F}_p \overset{d_p}\longrightarrow \mathbb{F}_{p-1} \cdots \overset{d_2}\longrightarrow \mathbb{F}_1 \overset{d_1}\longrightarrow \mathbb{F}_0 \overset{d_0}\longrightarrow M \longrightarrow 0 \qquad (*)$$ 
Here, $\mathbb{F}_i = \bigoplus_{\sigma \in \ZZ^n} S(-\sigma)^{\beta_{i,  \sigma}}$, where $S(-\sigma)$ denotes the free module obtained by shifting the degrees in $S$ by $\sigma$. The numbers $\beta_{i , \sigma}$'s are positive integers and are called the multigraded Betti numbers of $M$. We often identify $\sigma$ with the monomial whose exponent vector is $\sigma$. For example, over $K[x,y]$, we may write $\beta_{i, x^2y}(M)$ instead of $\beta_{i, (2,1)}(M)$.

For every $j \in \ZZ$, $\beta_{i,j}=\sum_{\{\sigma ~|~ |\sigma|=j\} } \beta_{i ,\sigma}$ is called the $(i,j)$-th standard graded Betti number of $M$. Three very important homological invariants that are related to these numbers are the Castelnuovo-Mumford regularity, or simply regularity, the depth and the projective dimension, denoted by $\text{reg}(M)$, $\text{depth}(M)$ and $\text{pd}(M)$ respectively:
\begin{align*}
\reg M & = \max \{|\sigma|-i ~|~ \beta_{i , \sigma} \neq 0 \} \\
\depth M & = \inf \{i ~|~\text{Ext}^i(K, M) \neq 0\} \\
\pd M & = \max \{i ~|~ \text{there is a } \sigma, \beta_{i , \sigma} \neq 0 \}.
\end{align*}

If $S$ is viewed as a standard graded $K$-algebra and $M$ is a graded $S$-module, then the graded Betti numbers of $M$ are also given by $\beta_{i,j}(M)=\dim_k \Tor_{i}(M,K)_j$, and so we have
\begin{align*}
\reg M & = \max \{j-i ~|~\Tor_{i} (M,K)_j \neq 0 \} \\
\pd M & = \max \{i~|~ \Tor_i (M,K) \neq 0 \} \\
\depth M & = n - \pd M.
\end{align*}

In practice, we often work with short exact sequences, so it is worthwhile to mention that the regularity can also be defined via the vanishings of \textit{local cohomology} modules with respect to the ``irrelevant maximal ideal''  $m=(x_1\dots ,x_n)$. Particularly,
for $i \geq 0,$ define
$$a_i(M):= \begin{cases}
\max  \{j~|~ H_m^i(M)_j \neq 0\} &  \text{if } H_m^i (M)\neq 0  \\
- \infty &  \text{otherwise.}
\end{cases}$$
Then, the regularity of $M$ is also given by $$\reg M = \max_i \{a_i(M)+i\}.$$
%  For a monomial ideal $I$ we say that $I$ is \emph{$k$-steps linear} whenever the minimal free resolution of $I$ over the polynomial ring
%is linear for $k$ steps, i.e., $\Tor_{i}^S(I,K)_j = 0$ for all $1\leq i\leq k$ and all $j\ne i+2s$. We say $I$ has linear minimal free resolution if the minimal free resolution is $k$-steps linear for all $k \geq 1$. We observe that an ideal is $k-$step linear when the matrices of the first $k-$ differentials of its linear minimal free resolution have linear entries.

\begin{example} \label{ex:resolution} Let $M=\displaystyle \mathbb{Q}[x_1, \dots ,x_5]/ (x_1x_2, x_2x_3, x_3x_4, x_4x_5, x_5x_1)$. Then the minimal free resolution of $M$ is:
$$0\longrightarrow \mathbb{F}_3 \overset{d_3}\longrightarrow \mathbb{F}_{2} \overset{d_2}\longrightarrow \mathbb{F}_1 \overset{d_1}\longrightarrow \mathbb{F}_0 \overset{d_0}\longrightarrow M \longrightarrow 0 $$ Here:
$$\beta_{0,\sigma}=1\text{ if }\sigma=1, \text{ and } \beta_{0, \sigma}=0 \text{ otherwise}$$
$$\beta_{1, \sigma}=1 \text{ if } \sigma= x_1x_2,x_2x_3,x_3x_4,x_4x_5,x_5x_1, \text{and } \beta_{1, \sigma}=0 \text{ otherwise}$$
$$\beta_{2 ,\sigma}=1 \text{ if } \sigma= x_1x_2x_3,x_2x_3x_4,x_1x_2x_5,x_1x_4x_5,x_3x_4x_5, \text{and } \beta_{2, \sigma}=0 \text{ otherwise}$$
$$\beta_{3, \sigma}=1\text{ if }\sigma=x_1x_2x_3x_4x_5, \text{ and } \beta_{3, \sigma}=0 \text{ otherwise}$$
\end{example}

An ideal $I$ in $S$ is said to be \textit{unmixed} if all its associated primes are minimal of same height. We say $I$ (or equivalently $S/I$)  is \textit{Cohen-Macaulay} if the Krull dimension and depth are equal. Cohen-Macaulay ideals are always unmixed and known to have many nice geometric properties. It can be checked that the module $M$ in Example \ref{ex:resolution} is Cohen-Macaulay.

\subsection{Algebraic objects with underlying combinatorial structures}
For any (hyper)graph $H$ over the vertex set $V(H) = \{x_1, \dots, x_n\}$, its \textit{edge ideal} is defined as follows:
$$I(H)=(\prod_{x \in e} x ~\big|~ e \in E(H)) \subseteq S = K[x_1, \dots, x_n].$$

\begin{example}
The edge ideal of the hypergraph, say $H,$ in Figure \ref{fig:hypergraph}, is
$$I(H)=(x_1x_2x_3x_4, x_3x_4x_5, x_5x_6x_7, x_7x_8,x_9).$$
\end{example}

For any simplicial complex $ \Delta$ over the vertex set $V = \{x_1, \dots, x_n\}$, its \textit{Stanley-Reisner ideal} is defined as follows:
$$I_{\Delta}=(\prod_{x \in e} x ~\big|~ e \text{ is a minimal non-face of } \Delta) \subseteq S.$$

\begin{example}
The Stanley-Reisner ideal  of the simplicial complex in Figure \ref{fig:simplicial} is
$$I_{\Delta}=(ac,ad,ae,bd,be,cf,df).$$
\end{example}

As we have seen before, (hyper)graphs are related to simplicial complexes via the notion of independence complex. The algebraic view of this relation exposes the following equality.

\begin{lemma}\label{lem.indcomplex}
Let $G$ be a  simple (hyper)graph and let $\Delta(G)$ be its independence complex. Then $$I(G)=I_{\Delta(G)}.$$
\end{lemma}

\begin{example}
The edge ideal of the graph $G$ in  Figure \ref{fig:independence} is the same as the the Stanley-Reisner ideal of the simplicial complex  $\Delta$ in Figure \ref{fig:simplicial}. Note that $\Delta=\Delta(G).$
$$I(G)= (ae,ad,bd,be,bf,ec,ef)$$
\end{example}

For any simple graph $G$ with $V(G)=\{x_1,\dots,x_n\}$, we define the $t$-\textit{path ideal} as: $$I_t(G)=(x_{i_1}\cdots x_{i_t} ~|~ i_k \neq i_l \text{ for } k \neq l,  x_{i_k}x_{i_{k+1}} \text{ an edge of } G).$$

\begin{example}
The $3$-path ideal of the five cycle $C_5: x_1x_2x_3x_4x_5$ is
$$I_3(C_5)=(x_1x_2x_3,x_2x_3x_4,x_3x_4x_5,x_4x_5x_1,x_5x_1x_2).$$
\end{example}
% An edge $e=uv$ in a graph $G$ is said to be a splitting edge if $I(G)=I(G_1) +(uv)$ is a splitting of $I(G)$, where $G_1!$ is the graph obtained from $G$ by deleting $e$.

Finally, as a matter of convention, for any (hyper)graph $G$, by $\reg G$ we mean $\reg I(G).$

\section{Formulas and inductive approaches} \label{chapter3}

Computing non-vanishings of local cohomology modules or Betti numbers of an ideal can be quite complicated.  As a result, inductive techniques that relate regularity of edge ideals with simplicial complexes and smaller ideals are employed as a common tool in the literature. In this context, induced structures such as induced subcomplexes and induced subgraphs have proven to be significant objects to investigate regularity of an edge ideal and its powers.

In this section, we focus on the methods that enable us to bound and compute regularity of an edge ideal and its powers. We start the section by recalling two important formulas. Then we address a few inductive bounds that are widely used in the literature.

\subsection{Hochster's and Takayama's Formulas}

Hochster's formula has been a significant tool in the study of squarefree monomial ideals due to its power to relate the  multigraded Betti numbers of a simplicial complex $\Delta$ to the non-vanishings of the reduced simplicial homology groups of $\Delta$ and its induced subcomplexes.

\begin{theorem}[\protect{Hochster's Formula \cite{Hochster}}] \label{formula.hochster}
Let $\Delta$ be a simplicial complex on the vertex set $V$ and let $I(\Delta)$ be its Stanley-Reisner ideal. Then

$$ \beta_{i,j} (I(\Delta))= \sum_{W \subseteq V, ~ |W| = j} \dim_k\big( \widetilde{H}_{j-i-2} (\Delta_W;K)\big)$$
where $\Delta_W$ is the restriction of $\Delta$ to the vertex set $W.$
\end{theorem}

For an arbitrary monomial ideal, Hochster's formula can not be employed. Takayama's formula given in \cite{Ta} can perform a similar task as Hochster's formula for this class of ideals.

Let $I \subseteq S = K[x_1, \dots, x_n]$ be a monomial ideal. Takayama's formula provides a combinatorial description for the non-vanishings of $\ZZ^n$-graded component $H_{\mm}^i (S/I)_{\bf a}$ for $\bf a \in \ZZ^n,$ and this description is given in terms of certain simplicial complexes $\Delta_{\bf a} (I)$ related to $I.$  Note that $S/I$ is an $\NN^n$-graded algebra, and so $H_{\mm}^i (S/I)$ is a $\ZZ^n$-graded module over $S/I.$

The simplicial complex $\Delta_{\bf a}(I)$ is called the \textit{degree complex} of $I.$ The construction of $\Delta_{\bf a}(I)$ was first given in \cite{Ta} and then simplified in  \cite{MT}. We recall the construction from \cite{MT}.

For $\textbf{a} =(a_1, \dots, a_n) \in \ZZ^n,$ set $x^a=x_1^{a_1} \cdots x_n^{a_n}$ and $ G_{\bf a} := \{ j \in \{1,\dots, n\} ~|~ a_j < 0\}.$ For every subset $ F \subseteq \{1,\dots, n\},$ let $S_F=S[x_j^{-1} ~|~ j \in F].$ Define

$$\Delta_{\bf a} (I)= \{ F \setminus G_{\bf a} ~|~ G_{\bf a} \subseteq F, x^{\bf a} \notin IS_F \}.$$

\begin{theorem}[\protect{Takayama's Formula \cite{Ta}}] \label{formula.takayama} Let $I \subseteq S$ be a monomial ideal, and let $\Delta(I)$ denote the simplicial complex corresponding to $\sqrt{I}$. Then

$$ \dim_k H_{\mm}^i (S/I)_{\bf a} =  \begin{cases}
      \dim_k\big( \widetilde{H}_{i-|G_{\bf a}|-1} (\Delta_{\bf a} (I),K\big) &  \textrm{ if } G_{\bf a} \in \Delta(I),\\
      0 & \textrm{ otherwise. }
   \end{cases} $$

\end{theorem}

The formula stated in \cite{MT} and here is different than the original formula introduced in \cite{Ta}.  The original formula has additional conditions on $\bf a$ for  $H_{\mm}^i (S/I)_{\bf a} =0.$ It follows from the proof in \cite{Ta} that those conditions can be omitted.

Due to the equality between edge ideals of graphs and Stanley-Reisner ideals of their independence complexes, we can employ Hochster's and Takayama's formulas in the study of regularity of edge ideals.

In order to deal with arbitrary monomial ideals, polarization is proved to be a powerful process to obtain a squarefree monomial ideal from a given monomial ideal.   For details of polarization we refer the reader to \cite{HH}.
\begin{definition}\label{pol_def}
Let $M=x_1^{a_1}\dots x_n^{a_n}$ be a monomial in  $S=K[x_1,\dots,x_n]$. Then we define the squarefree monomial $P(M)$ ({\it polarization} of $M$) as
$$P(M)=x_{11}\dots x_{1a_1}x_{21}\dots x_{2a_2}\dots x_{n1}\dots x_{na_n}$$ in the polynomial ring $R=K[x_{ij} \mid 1\leq i\leq n,1\leq j\leq a_i]$.
If $I=(M_1,\dots,M_q)$ is an ideal in $S$, then the polarization of $I$, denoted by $I^{\pol}$, is defined as $I^{\pol}=(P(M_1),\dots,P(M_q))$.
\end{definition}

The regularity is preserved under polarization.

\begin{corollary} \cite[Corollary 1.6.3.d]{HH}  Let $I \subset S$ be a monomial ideal and $I^{\pol}\subset R$ be its polarization. Then $\reg (S/I) =\reg (R/I^{\pol}).$
\end{corollary}

%%%%%%%%%%%%%%%%%%%%%%%%%%%%%%%%%%%%%%%%

\subsection{Inductive techniques}

We start the section with an easy yet essential consequence of Hochster's formula that links the regularity of a graph with regularity of its induced subgraphs.

\begin{lemma}\label{lem.induced}
Let $G$ be a simple graph. Then $\reg H \leq \reg G$ for any induced subgraph $H$ of $G.$
\end{lemma}

For any homogeneous ideal $I \subseteq S$ and any homogeneous element $m$ of degree $d,$ the following short exact sequences are used as standard tools in commutative algebra and also proved to be very useful in computing regularity of edge ideals and their powers.

\begin{eqnarray}\label{eq.colonadd}
0 \longrightarrow \frac{R}{I:m} (-d) \xrightarrow{\cdot m} \frac{R}{I} \longrightarrow \frac{R}{I+(m)} \longrightarrow 0
\end{eqnarray}

Let  $I$ and $J$ be ideals in $S.$ Another useful exact sequence is:

\begin{eqnarray}\label{eq.intersectadd}
0 \longrightarrow \frac{R}{I \cap J}  \rightarrow \frac{R}{I}  \oplus \frac{R}{J} \longrightarrow \frac{R}{I+J} \longrightarrow 0
\end{eqnarray}

We can see how regularity changes in a short exact sequence by taking the associated long exact sequence of local cohomology modules. Particularly,  we have the following inductive bound (the second statement is the content of \cite[Lemma 2.10 ]{DHS}).

\begin{lemma}\label{lem.regineq}
Let $I \subseteq S$ be a monomial ideal, and let $m$ be a monomial of degree $d$.  Then
\[
\reg I  \leq \max\{ \reg (I : m) + d, \reg (I,m)\}.
\]
Furthermore, if $x$ is a variable appearing in $I,$ then
$$\reg I \in \{\reg (I:x) + 1, \reg (I,x)\}.$$
\end{lemma}

In the statement of Lemma \ref{lem.regineq}, the expression $x$ is a variable appearing in $I$ means some of the minimal generators of $I$ is divisible by $x.$ Note that if $x$ is a variable not appearing in $I,$ then $\reg (I,x)=\reg I.$ In case of edge ideals, if $x$ is an isolated vertex in $G,$ we can drop the vertex $x$ when computing regularity. Thus $\reg (I(G):x) = \reg I(G \setminus N_G[x])$ and $\reg (I(G),x) = \reg (I(G \setminus x))$ for any vertex $x$ in $G.$  Then Lemma \ref{lem.regineq} can be restated in terms of edge ideals.

\begin{lemma}\label{lem.deletecontract}
Let $x$ be a vertex in $G.$ Then
 $$\reg G \in \{ \reg (G \setminus N_G[x]) +1, \reg (G \setminus x)\}.$$
\end{lemma}

Kalai and Meshulam  \cite{KM} proved the following result for squarefree monomial ideals and Herzog \cite{Herzog} generalized it to any monomial ideal.

\begin{theorem} \label{thm.kalaimeshulam}
Let $I_1, \dots, I_s$ be squarefree monomial ideals in $S.$ Then

$$ \reg \Big( S \big/ \sum_{i=1}^s I_i \Big) \leq \sum_{i=1}^s \reg S /I_i$$

\end{theorem}

In the case of edge ideals, we have the following bound.

\begin{corollary}\label{cor.kalaimeshulam}
Let $G$ be a simple graph. If $G_1, \dots, G_s$ are subgraphs of $G$ such that $E(G)= \bigcup_{i=1}^s E(G_i)$ then
$$ \reg \Big( S /  I(G) \Big) \leq \sum_{i=1}^s \reg S /I(G_i)$$

\end{corollary}

If $G$ is the disjoint union of graphs $G_1, \dots, G_s$ then the above equality is achieved by using K\"unneth formula in algebraic topology.

\begin{corollary}\label{cor.disjoint}
Let $G$ be a simple graph. If $G$ can be written as a disjoint union of graphs $G_1, \dots, G_s$  then
$$ \reg \Big( S /  I(G) \Big) = \sum_{i=1}^s \reg S /I(G_i)$$

\end{corollary}

In the study of powers of edge ideals, Banerjee developed the notion of even-connection and gave an important inductive inequality in \cite{Ba}.

\begin{theorem} \label{thm.Banerjee}
For any simple graph $G$ and any $s \geq 1,$ let the set of minimal monomial generators of $I(G)^s$ be $\{m_1, \dots, m_k\}.$ Then

$$ \reg I(G)^{s+1} \leq \max \{\reg (I(G)^{s+1}: m_l)+2s , 1 \leq l \leq k, \reg I(G)^s \}.$$

\par In particular, if for all $s \geq 1$ and for all minimal monomial generators $m$ of $I(G)^s,$ $\textrm{reg }(I(G)^{s+1}: m) \leq 2$ and $\reg I(G) \leq 4,$ then $\reg I(G)^{s+1}= 2s+2$ for all $s \geq 1;$  as a consequence $I(G)^{s+1}$ has a linear minimal free resolution.
\end{theorem}

Working with the above inequality requires a good understanding of the ideal $(I(G)^{s+1}: m)$ and its generators when $m$ is a minimal monomial generator of $I(G)^s.$  The even-connection definition is key to attain this goal. We recall the definition of even-connectedness and its important properties from \cite{Ba}.

\begin{definition}\label{even_connected} Let $G=(V,E)$ be a graph. Two vertices $u$ and $v$ ($u$ may be the same as $v$) are
said to be even-connected with respect to an $s$-fold product $e_1\cdots e_s$ where $e_i$'s are edges of $G$, not necessarily distinct,
if there is a path $p_0p_1\cdots p_{2k+1}$, $k\geq 1$ in $G$ such that:
\begin{enumerate}
 \item $p_0=u,p_{2k+1}=v.$
 \item For all $0 \leq l \leq k-1,$ $p_{2l+1}p_{2l+2}=e_i$ for some $i$.
 \item For all $i$, $ \mid\{l \geq 0 \mid p_{2l+1}p_{2l+2}=e_i \}\mid
   ~ \leq  ~ \mid \{j \mid e_j=e_i\} \mid$.
 \item For all $0 \leq r \leq 2k$, $p_rp_{r+1}$ is an edge in $G$.
\end{enumerate}
\end{definition}

Fortunately, it turns out that $(I(G)^{s+1} : m)$ is generated  by monomials in degree 2.

\begin{theorem}\label{even_connec_equivalent}\cite[Theorem 6.1 and Theorem 6.7]{Ba} Let $G$ be a graph with edge ideal
$I = I(G)$, and let $s \geq 1$ be an integer. Let $m$ be a minimal generator of $I^s$.
Then $(I^{s+1} : m)$ is minimally generated by monomials of degree 2, and $uv$ ($u$ and $v$ may
be the same) is a minimal generator of $(I^{s+1} : m )$ if and only if either $\{u, v\} \in E(G) $ or $u$ and $v$ are even-connected with respect to $m$.
 \end{theorem}

After polarization, the ideal $(I(G)^{s+1} : m)$ can be viewed as the edge ideal of a graph that is obtained from $G$ by adding even-connected edges with respect to $m.$ Some of the combinatorial properties of this ideal in relation to $G$ are studied in \cite{JNS}.

\section{Regularity of edge ideals} \label{chapter4}

Computing and/or bounding the regularity of edge ideals is the foundation in the study of regularity of powers of edge ideals. In this section, we identify the combinatorial structures that are related to regularity of edge ideals. First, we collect the general upper and lower bounds given for the regularity of edge ideals, then we present the list of known classes of graphs where the regularity is computed explicitly.

Note that regularity of an edge ideal is bounded below by 2, which is the generating degree of an edge ideal. Thus, identifying combinatorial structures of a graph with regularity 2 can be considered as the base case of results in this section. The following combinatorial characterization of such graphs is nowadays often referred to as \emph{Fr\"oberg's characterization}. It was, in fact, given first in topological language by Wegner \cite{Wegner} and later, independently, by Lyubeznik \cite{L} and Fr\"oberg  \cite{Fr} in monomial ideals language.

\begin{theorem}[\protect{\cite[Theorem 1]{{Fr}}}] \label{thm.regtwo}
Let $G$ be a simple graph. Then $\reg I(G) = 2$ if and only if $G$ is a co-chordal graph.
\end{theorem}

\subsection{Lower and upper bounds}

One of the graph-theoretical invariants that can be related to the regularity is the induced matching number.  The first result revealing this relation is due to Katzman and it provides a general lower bound on the regularity of edge ideals.

\begin{theorem}[\protect{\cite[Lemma 2.2]{Katzman}}] \label{thm.reggenerallower}
Let $G$ be a simple graph and $\nu(G)$ be the maximum size of an induced matching in $G.$ Then
$$\reg I(G) \ge \nu(G)+1.$$
\end{theorem}

\begin{sketch}
Let $\{e_1, \dots, e_r\}$ be an induced matching of maximal size in $G.$ Suppose $H$ is the induced subgraph of $G$ with $E(H)=\{e_1, \dots, e_r\}.$ Note that all the edges in $H$  are disjoint. Thus $\reg (H) =r+1.$ By Lemma \ref{lem.induced}, $\reg (G) \ge \nu(G)+1.$
\end{sketch}

Another graph-theoretical invariant of interest is the matching number and this invariant actually emerges as a general upper bound for any graph.

\begin{theorem}[\protect{\cite[Theorem 6.7]{HVT2008}, \cite[Theorem 11]{Russ}}] \label{thm.regupper}
Let $G$ be a simple graph. Let $\beta(G)$ be the minimum size of a maximal matching in $G.$ Then
$$\reg I(G) \le \beta(G)+1.$$
\end{theorem}

\begin{example}
Let $G$ be the graph given in Figure \ref{fig:maxmatch}. It is clear that $\nu(G)=\beta(G)=3.$ Then $\reg G=4$ by Theorem \ref{thm.reggenerallower} and  Theorem \ref{thm.regupper}. In general, if $G$ is a simple graph with disjoint edges, then the above bounds coincide and become an equality.

\begin{figure}[h]
  \includegraphics[width=0.26\linewidth]{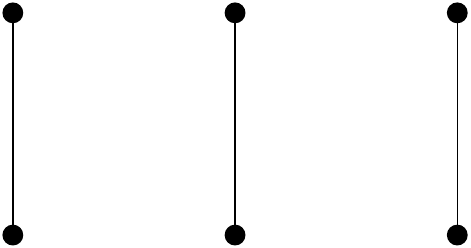}
  \caption{A graph with disjoint edges}
  \label{fig:maxmatch}
\end{figure}
\end{example}

Comparing the existing lower and upper bounds yields to interesting classifications, particularly when these bounds coincide. For example, Cameron and Walker \cite[Theorem 1]{CW} gave the first classification of graphs $G$ with $\nu(G)=\beta(G).$ Then Hibi et al  \cite{HHKO} modified their result slightly and gave a full generalization with some corrections.

\begin{example}
Let $G$ be the graph given in Figure \ref{fig:cochordal}.  One can easily verify that $ \beta(G)=2.$ However $\reg G= 2$ by Fr\"oberg's characterization.  In general, we can find examples of graphs where $\beta (G)+1$ can be arbitrarily large compared to $\reg G.$
\begin{figure}[h]
  \includegraphics[width=0.5\linewidth]{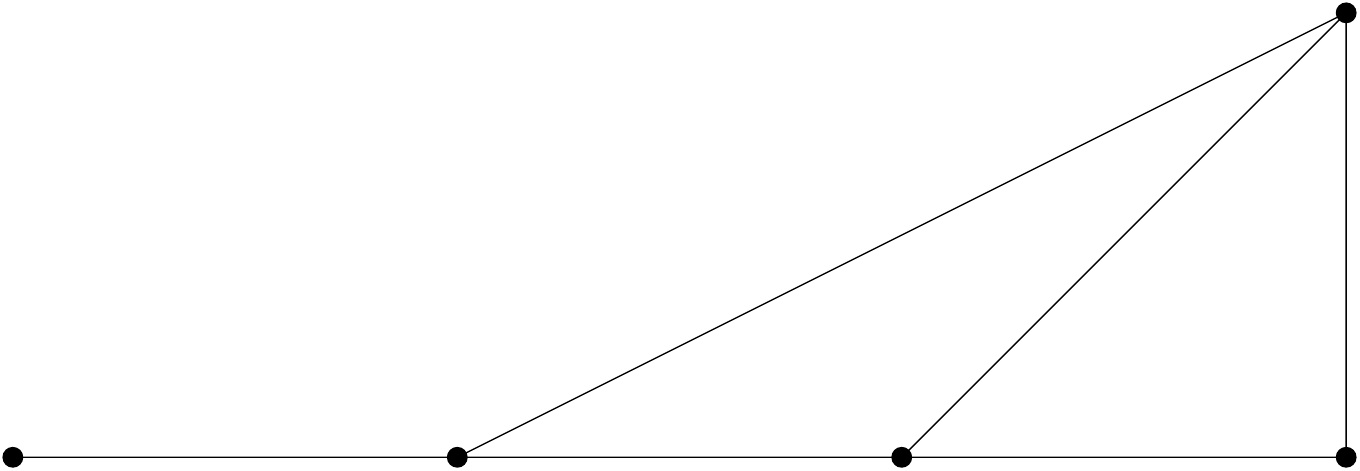}
  \caption{A co-chordal graph}
  \label{fig:cochordal}
\end{figure}
\end{example}

The above upper bound is strengthened by making use of co-chordal subgraph covers of a graph and this bound is proved by Woodroofe. Recall that a graph is co-chordal if its complement is chordal and $\text{co-chord}(G)$, the \emph{co-chordal number} \index{co-chordal number} of $G$,
denotes the least number of co-chordal subgraphs of $G$ whose union is $G$. The graph in Figure \ref{fig:cochordal} is an example of a co-chordal graph.

Let $\{e_1, \dots, e_r\}$ be a maximal matching of minimal size in $G.$  For each $i,$ let $G_i$ be the subgraph of $G$ with edges $e_i \cup \{  \text{ edges in } G \text{ adjacent to } e_i\}.$ Note that $G_1, \dots, G_r$ forms a co-chordal subgraph cover of $G,$ thus $\text{co-chord}(G) \le \beta (G)$ and the bound in Theorem \ref{thm.reggeneralupper} improves the bound in Theorem \ref{thm.regupper}.

\begin{theorem}[\protect{\cite[Lemma 1]{Russ}}] \label{thm.reggeneralupper}
\index{regularity}
\index{co-chordal number}
Let $G$ be a simple graph. Then
$$\reg I(G) \le \text{co-chord}(G)+1.$$
\end{theorem}

\begin{sketch}
 Let  $\text{co-chord}(G)=r$ and let $G_1, \dots, G_r$ be a co-chordal cover of $G$. It follows from Fr\"oberg's characterization of regularity 2 graphs, Theorem \ref{thm.regtwo}, that $\reg (G_i) =2$ for each $i\in \{1, \dots, r\}.$ Then the result follows immediately from Corollary \ref{cor.kalaimeshulam}.
\end{sketch}

Gap-free graphs have also been of interest in the investigation of regularity. These graphs arise naturally since their induced matching number is 1 which is the smallest possible number it could be. However, computing or bounding the regularity  of this class of graphs is not so easy.  Furthermore, there are very few examples on how large the regularity of such edge ideals can be, see  \cite{NP} by Nevo and Peeva for an example of a gap-free graph $G$ in 12 variables with $\reg G= 4.$ Putting an additional condition on the gap-freeness of $G$ may result with an upper bound and  it is indeed achieved in \cite{Nevo} .

\begin{theorem}[\protect{\cite[Theorem 1.2]{Nevo}, \cite[Proposition 19]{DHS}}] \label{thm.gapclawfree}
If  $G$ is gap-free and claw-free, then  $\reg I(G) \le 3.$
\end{theorem}

\begin{sketch} Let $x$ be a vertex in $G$ with the maximum possible degree. By Lemma \ref{lem.deletecontract},  we have $\reg (G) \leq \max \{ \reg (G \setminus N_G[x]) +1, \reg (G \setminus x)\}.$ Note that induced subgraphs $G \setminus N_G[x]$ and $G\setminus x$ of $G$ are both gap-free and claw-free. It follows from the induction on the number of vertices that $\reg (G\setminus x) \le 3.$ Thus it suffices to show that $\reg ( G\setminus N_G[x]) \le 2.$ By  Fr\"oberg's characterization, it is enough to prove that $( G\setminus N_G[x])^c$ is chordal and it is proved by contradiction:
 \begin{enumerate} 
 \item Suppose  $(G\setminus N_G[x])^c$ has an induced cycle on $w_1,w_2, \ldots, w_n$ (in order) of length at least 4. 
\item Any vertex of $G$ is distance 2 from $x$ in $G$ by  \cite[Proposition 3.3]{DHS}, then $\{x,y\}$ and $\{y,w_1\}$ are edges in $G$ for some vertex $y.$ 
\item Note that either $\{y,w_2\}$ or $ \{y,w_n\}$ must be an edge in $G.$ Otherwise edges $ \{w_2,w_n\}$ and $\{x,y\}$ form a gap in $G.$ 
\item  Without loss of generality, suppose $\{y,w_2\}$ is an edge. Then the induced subgraph on $\{x,y,w_1,w_2\}$ is a claw in $G,$ a contradiction.
 \end{enumerate}
\end{sketch}

\begin{example}
There are examples of gap-free and claw-free graphs where both values in  Theorem  \ref{thm.gapclawfree} can be attained for regularity. For example, if $G^c$ is a tree then $\reg G=2$ and if $G^c$ is $C_n$ for $n\geq 5$ then $\reg G= 3.$
\end{example}

This result is generalized to $n$-claw free and gap-free graphs by Banerjee in \cite{Ba}. The proof of the general case follows similarly and uses induction on $n.$

\begin{theorem}[\protect{\cite[Theorem 3.5]{Ba}}] \label{thm.gapnclawfree}
If  $G$ is gap-free and $n$-claw-free, then  $\reg I(G) \le n.$
\end{theorem}

Another special class of graphs are planar graphs. This class of graphs emerges frequently in applications since they can be drawn in the plane without edges crossing. It is given in \cite{Russ} that even though regularity of a planar graph may be arbitrarily large, regularity of its complement can be bounded above by 4.

\begin{theorem}[\protect{\cite[Theorem 3.4]{Russ}}] \label{thm.regplanarupper}
If  $G$ is a planar graph, then  $\reg I(G^c) \le 4.$
\end{theorem}

%%%%%%%%%%%%%%%%%%%%%%%%%%%%%%%%%%%%%%%%

\subsection{Exact values}

Computing the regularity for special classes of graphs has been an attractive research topic in the recent years. To that extend characterizing edge ideals with certain regularity has been of interest as well. However, very little is known in the latter case.

A combinatorial characterization of edge ideals with regularity 3 is still not known. However, a partial result for bipartite graphs is achieved by Fern$\acute{\textrm{a}}$ndez-Ramos and Gimenez in \cite{FG}. Recall that a graph $G$ is bipartite if the vertices $V$ can be partitioned into disjoint subsets $V=  X \cup Y$ such that $\{x,y\} $ is an edge in $G$ only if $x \in X$ and $y \in Y$ or vice versa. The \textit{bipartite complement} of a bipartite graph $G,$ denoted by $G^{bc},$ is the bipartite graph over the same partition of vertices and $ \{x,y\} \in G^{bc}$ if and only if  $\{x,y\} \notin G$ for   $x \in X, y \in Y.$

\begin{theorem}[\protect{\cite[Theorem 3.1]{FG}}] \label{thm.regbipartitethree}
Let $G$ be a connected bipartite graph. Then $\reg I(G)=3$ if and only if $G^c$ has no induced cycles of length $\ge 4$ and $G^{bc}$ has no induced cycles of length $\ge 6.$
\end{theorem}

In most of the known cases, it turns out that regularity can be expressed in terms of the induced matching number. We first recall the results when the regularity is one more than the induced matching number.

\begin{theorem}\label{thm.inducedmatch1}
Let $G$ be a simple graph and $\nu(G)$ be the induced matching number of $G.$ Then

$$\reg I(G) = \nu(G)+1$$

in the following cases:

\begin{enumerate}[label=(\alph*)]
\item $G$ is a chordal graph (see \cite{HVT2008});
\item $G$ is a weakly chordal graph (see \cite{Russ});
\item $G$ is sequentially Cohen-Macaulay bipartite graph (see  \cite{Van});
\item $G$ is unmixed bipartite graph  (see \cite{Ku});
\item $G$ is very well-covered graph  (see \cite{MMCRTY});
\item $G$ is vertex decomposable graph and has no closed circuit of length 5 (see \cite{KMo});
\item $G$ is a $(C_4,C_5)$-free vertex decomposable graph (see \cite{BC});
\item $G$ is a unicyclic graph with cycle $C_n$ when
\begin{enumerate}[label=(\roman*)]
\item $ n \equiv 0, 1 ~( mod ~3)$ (see \cite{ ABS, BHT, J}) or
\item $ \nu(G \setminus \Gamma(G)) < \nu (G)$ where $\Gamma(G)$ is the collection of all neighbors of the roots in the rooted trees attached to $C_n $  (see \cite{ABS}).
\end{enumerate}
\end{enumerate}
\end{theorem}

\begin{figure}[h]
  \includegraphics[width=0.33\linewidth]{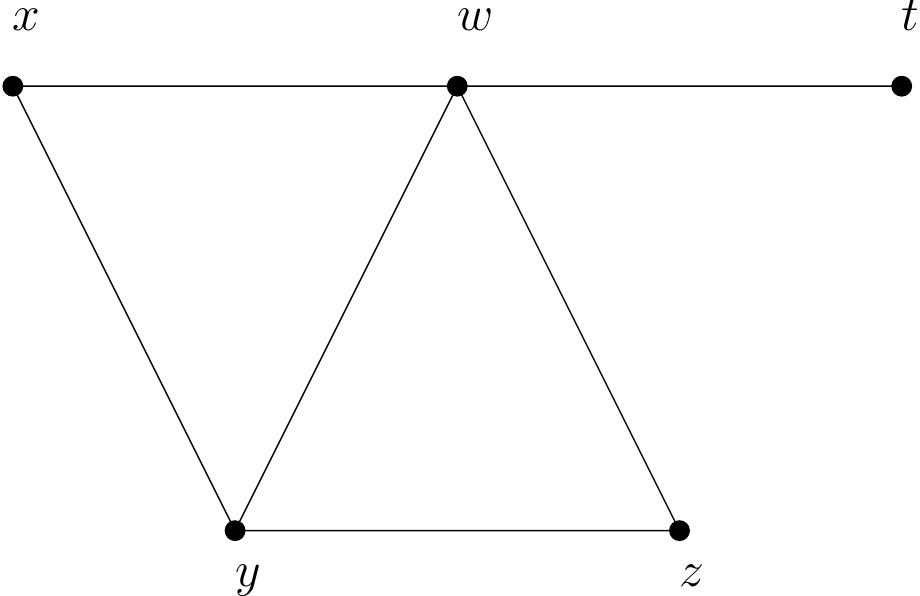}
  \caption{A vertex decomposable and a (weakly) chordal graph}
  \label{fig:vertexdec}
\end{figure}

\begin{remark} Chordal graphs are vertex decomposable and vertex decomposable graphs are sequentially Cohen-Macaulay (see \cite{FVT, Russ2}).
Sequentially Cohen-Macaulay bipartite graphs are vertex decomposable (see \cite{Van}). Also if $G$ is a very well-covered graph then $G$ is unmixed and if $G$ is chordal, then it is weakly chordal. Thus, $(b)$ implies $(a)$, and $(e)$ implies $(d)$. Note that bipartite graphs have no odd cycles. Hence, $(f)$ implies $(c)$.
\end{remark}

\begin{example}
Let $G$ be the graph given in Figure \ref{fig:unicyclic1}. Then $\Gamma(G)=\{x_6, x_7,x_8\}$ and $\nu(G \setminus \Gamma(G))=2 < \nu(G)=3.$
Furthermore, this graph does not belong to any of the classes described in $(a)$---$(g)$.
\begin{figure}[h]
  \includegraphics[width=0.6\linewidth]{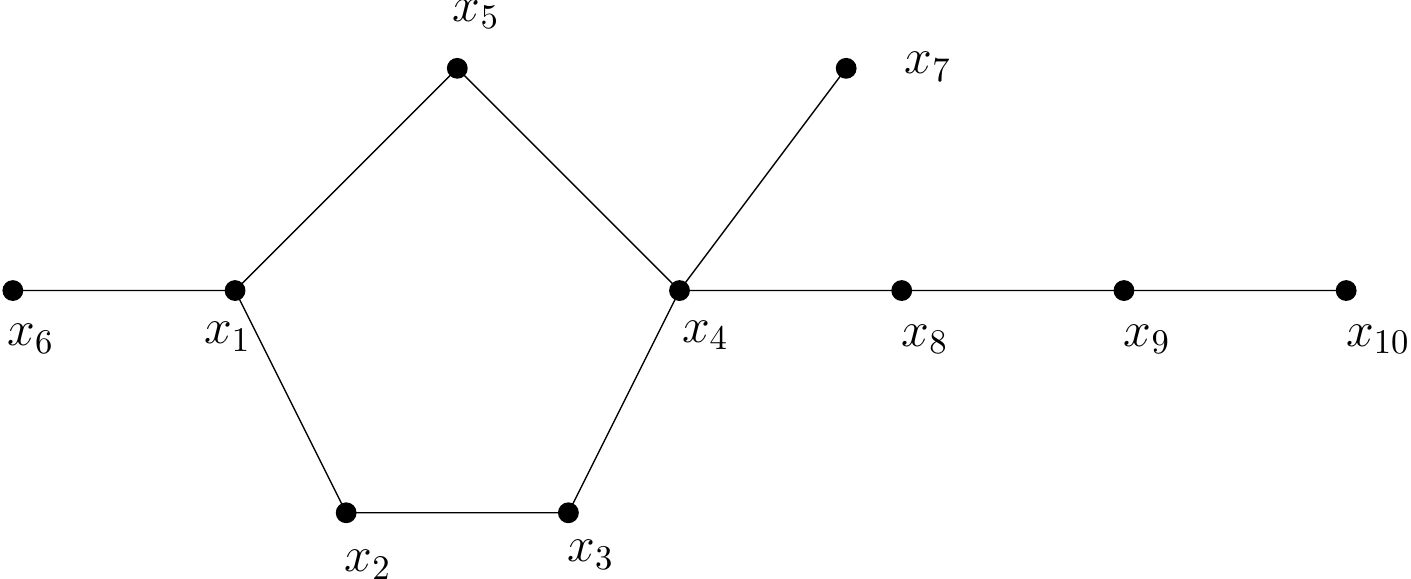}
  \caption{A graph satisfying (h) (ii) in Theorem \ref{thm.inducedmatch1}}
  \label{fig:unicyclic1}
\end{figure}
\end{example}

It is of interest to find different expressions for the regularity of edge ideals.  The next  result collects all known cases in which the regularity has a different expression than the above classes and it is still in terms of induced matching.

\begin{theorem} \label{thm.inducedmatch2}
Let $G$ be a simple graph and $\nu(G)$ be the induced matching number of $G.$ Then

$$\reg I(G) = \nu(G)+2$$

in the following cases:

\begin{enumerate}[label=(\alph*)]
\item $G$ is an $n$-cycle $C_n$ when $ n \equiv 2 ~ (mod ~3)$ (see \cite{ BHT, J});
\item $G$ is a unicyclic graph with cycle $C_n$ when  $ n \equiv 2 ~( mod ~3)$ and $ \nu(G \setminus \Gamma(G)) = \nu (G)$ where $\Gamma(G)$ is the collection of all neighbors of the roots in the rooted trees attached to $C_n $  (see \cite{ABS}).
\end{enumerate}
\end{theorem}

\section{Regularity of powers of edge ideals} \label{chapter5}

The regularity of powers of an edge ideal is considerably harder to compute than that of the edge ideal itself. However, in many cases, known bounds and exact formulas for the regularity of $I(G)$ inspire new bounds and exact formulas for the regularity of $I(G)^q$, for $q \ge 1$. This section is divided into two parts, where in the first subsection we list a number of lower and upper bounds, and in the second subsection we give the exact values of $\reg I(G)^q$, for $q \ge 1$, for special classes of graphs.

\subsection{Lower and upper bounds} Just like with studies on the regularity of edge ideals, the induced matching number of a graph is ultimately connected to the regularity of powers of its edge ideal. The following general lower bound generalizes that of Theorem \ref{thm.reggenerallower} for an edge ideal to its powers.

\begin{theorem}[\protect{\cite[Theorem 4.5]{BHT}}] \label{thm.generallower}
Let $G$ be any graph and let $I = I(G)$ be its edge ideal.
Then for all $q \ge 1$, we have
$$\reg I^q \ge 2q+\nu(G)-1.$$
\end{theorem}

\begin{sketch} The proof is based on the following observations:
\begin{enumerate}
\item If $H$ is an induced subgraph of $G$ then for any $i, j \in \ZZ$, we have
\begin{align}
\beta_{i,j}(I(H)^q) \le \beta_{i,j}(I(G)^q). \label{eq.500}
\end{align}
In particular, this gives $\reg I(H)^q \le \reg I(G)^q$.
\item If $H$ is the induced subgraph of $G$ consisting of a maximal induced matching of $\nu(G)$ edges then for any $q \ge 1$, we have
\begin{align}
\reg I(H)^q = 2q+\nu(G)-1. \label{eq.501}
\end{align}
\end{enumerate}

To establish the first observation, recall that the \emph{upper-Koszul simplicial complex} $K^\alpha(I)$ associated to a monomial ideal $I \subseteq S = K[x_1, \dots, x_n]$ at degree $\alpha \in \ZZ^n$ consists of faces
$$\Big\{W \subseteq \{1, \dots, n\} ~\Big|~ \dfrac{x^\alpha}{\prod_{w \in W} w} \in I\Big\},$$
and a variation of Hochster's formula (Theorem \ref{formula.hochster}) gives
$$\beta_{i,\alpha}(I) = \dim_k \widetilde{H}_{i-1}(K^\alpha(I);K) \text{ for all } i \ge 0.$$
The inequality (\ref{eq.500}) then follows by noting that for $\alpha \in \ZZ^n$ with $\supp(\alpha) \subseteq V_H$, $K^\alpha(I(H)^q) = K^\alpha(I(G)^q).$

The second observation is proved by induction, noting that in this case $I(H)$ is a complete intersection.
\end{sketch}

\begin{figure}[h]
  \includegraphics[width=0.4\linewidth]{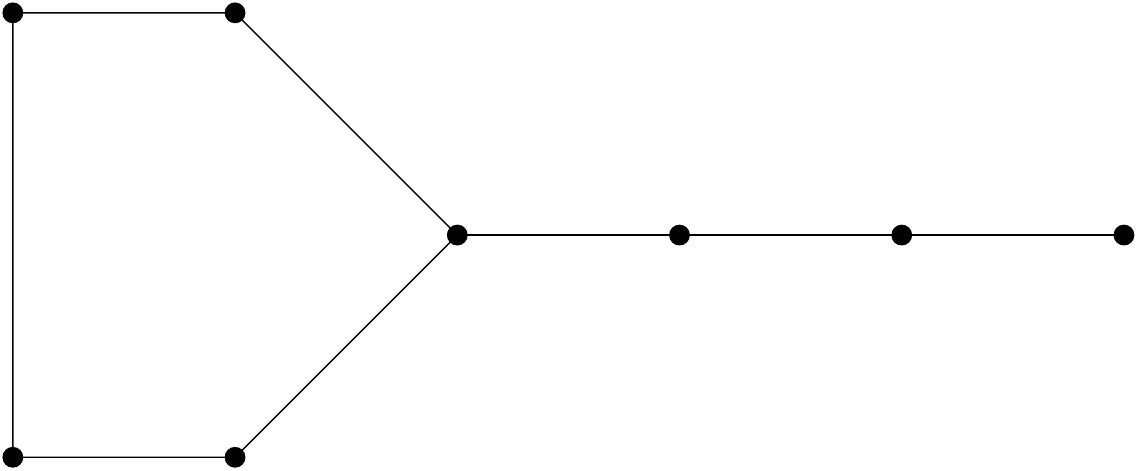}
  \caption{A uniyclic graph}
  \label{fig:lowerbound}
\end{figure}

\begin{example}
Though there are many classes of graphs for which the lower bound given in Theorem \ref{thm.generallower} is attained, there are also classes of graphs where the asymptotic linear function $\reg I(G)^q$ is strictly bigger than $2q + \nu(G)-1$ for all $q \gg 0$. Let $G$ be the graph depicted in Figure \ref{fig:lowerbound}. Then it is easy to see that $\nu(G)=2$, whereas $\reg I(G)^q= 2q+2$ for $q\ge 1.$
\end{example}

A similar general upper bound generalizing that of Theorem \ref{thm.reggeneralupper} is, unfortunately, not available. It is established, by Jayanthan, Narayanan and Selvaraja \cite{JNS}, only for a special class of graphs --- bipartite graphs.

\begin{theorem}[\protect{\cite[Theorem 1.1]{JNS}}] \label{thm.JNS}
Let $G$ be a bipartite graph and let $I = I(G)$ be its edge ideal.
Then for all $q \ge 1$, we have
$$\reg I^q \le 2q+\text{co-chord}(G)-1.$$
\end{theorem}

\begin{sketch} The statement is proved by induction utilizing Theorem \ref{thm.Banerjee}. The crucial step in the proof is to show that for any $q \ge 1$ and any collection of edges $e_1, \dots, e_q$ of $G$ (not necessarily distinct), we have
\begin{align}
\reg (I^{q+1}:e_1\dots e_q) \le \text{co-chord}(G)+1. \label{eq.503}
\end{align}

Observe that when $G$ is a bipartite graph, by \cite[Lemma 3.7]{AB},
$$I^{q+1}:e_1 \dots e_q = (((I^2:e_1)^2: \dots)^2 : e_q),$$
and so (\ref{eq.503}) itself can be obtained by induction. To this end, let $G'$ be the graph associated to the polarization of $I^2:e$ for an edge $e$ in $G$ (which is generated in degree 2 by Theorem \ref{even_connec_equivalent}). The proof is completed by establishing the following facts:
\begin{enumerate}
\item $\reg I(G') \le \text{co-chord}(G')+1.$ This inequality was proved to hold for any graph in \cite{Russ}.
\item $\text{co-chord}(G') \le \text{co-chord}(G).$ This is combinatorial statement, which can be shown by analyzing how $G'$ is constructed from $G$.
\end{enumerate}
\end{sketch}

\begin{remark}
Let $G$ be a 4-cycle. We can easily verify that $\nu(G)=1$ and $\text{co-chord}(G)=1.$ Since the lower bound in Theorem \ref{thm.generallower} coincides with the upper bound in Theorem \ref{thm.JNS}, we have $\reg I(G)^q=2q= 2q+\nu(G)-1=2q+ \text{co-chord}(G)-1$ for all $q \ge 1$.
Classes of graphs for which the two upper and lower bounds agree were discussed in \cite[Corollary 5.1]{JNS}.

On the other hand, the upper bound given in Theorem \ref{thm.JNS} can be strict. For example, if $G$ is $C_8,$ then $\nu(G)=2$ and $\text{co-chord}(G)=3.$ By \cite[Theorem 5.2]{BHT}, it is known that $\reg I(G)^q=2q+1 < 2q + \text{co-chord}(G)-1$ for all $q \ge 1$.
\end{remark}

For a special class of graphs --- gap-free graphs --- there is another upper bound that was proved by Banerjee \cite{Ba}.

\begin{theorem}[\protect{\cite[Theorem 6.19]{Ba}}] \label{thm.Ba_gapfree}
Let $G$ be a gap-free graph and let $I = I(G)$ be its edge ideal. Then for all $q \ge 2$, we have
$$\reg I^q \le 2q + \reg I - 1.$$
\end{theorem}

\begin{remark} The bound in Theorem \ref{thm.Ba_gapfree} is slightly weaker than the conjectural bound of Conjecture \ref{conj.BBH}.
\end{remark}

%%%%%%%%%%%%%%%%%%%%%%%%%%%%%%%%%%%%%%%%

\subsection{Exact values} In most of the known cases where the asymptotic linear function $\reg I(G)^q$ can be computed explicitly, the lower bound in Theorem \ref{thm.generallower} turns out to give the exact formula. In this subsection, we describe those instances.

In \cite[Theorem 3.2]{HHZ}, the authors showed that $I(G)$ has a linear resolution if and only if $I(G)^q$ has a linear resolution for all $q\geq 1.$ Combination of their result with Fr\"oberg's characterization yields to the exact value of regularity of powers of co-chordal graphs. 

\begin{theorem} \label{thm.linear}
Let $G$ be a co-chordal graph and $I = I(G)$ be its edge ideal. Then for all $q \ge 1$, we have
$$\reg I^q = 2q.$$
\end{theorem}

As often the case when dealing with graphs, one of the first classes to consider is that of trees and forests.

\begin{theorem}[\protect{\cite[Theorem 4.7]{BHT}}] \label{thm.forest}
Let $G$ be a forest and let $I = I(G)$ be its edge ideal.
Then for all $q \ge 1$, we have
$$\reg I^q = 2q+\nu(G) -1.$$
\end{theorem}

\begin{sketch} Thanks to the general lower bound in Theorem \ref{thm.generallower}, it remains to establish the upper bound
$$\reg I^q \le 2q+\nu(G)-1.$$
In fact, a more general inequality can be proved for induced subgraphs of $G.$ Let $H$ and $K$ be induced subgraphs of $K$ such that
$$E_H \cup E_K = E_G \text{ and } E_H \cap E_K = \emptyset.$$
Then the required upper bound follows from the following inequality (by setting $H$ to be the empty graph):
\begin{align}
\reg (I(H) + I(K)^q) \le 2q + \nu(G)-1. \label{eq.504}
\end{align}

The inequality (\ref{eq.504}) is proved by induction on $q+|V_K|$, making use of the short exact sequence arising by taking quotient and colon with respect to a leaf $xy$ in $K$. Note that, by \cite[Lemma 2.10]{MO}, we have
$$(I(H)+I(K)^q) : xy = (I(H):xy) + (I(K)^q:xy) = (I(H):xy) + I(K)^{q-1}.$$

\end{sketch}

The next natural class of graphs to consider is that of cycles and graphs containing exactly one cycle (i.e., \emph{unicyclic} graphs).

\begin{theorem}[\protect{\cite[Theorem 5.2]{BHT}}] \label{thm.cycle}
Let $C_n$ denote the $n$-cycle and let $I = I(C_n)$ be its edge ideal.
Let $\nu = \lfloor \frac{n}{3} \rfloor$ be the induced matching number
of $C_n$. Then
$$\reg I = \left\{ \begin{array}{lcll} \nu+1 & \text{if} & n \equiv 0,1 & (\text{mod } 3) \\
\nu+2 & \text{if} & n \equiv 2 & (\text{mod } 3) \end{array} \right.$$
and for any $q \ge 2$, we have
$$\reg I^q = 2q+\nu-1.$$
\end{theorem}

\begin{sketch} The first statement was already proved in \cite{J}. To prove the second statement, again thanks to the general lower bound of Theorem \ref{thm.generallower}, it remains to establish the upper bound
\begin{align}
\reg I^q \le 2q+\nu-1. \label{eq.505}
\end{align}

The inductive method of Theorem \ref{thm.Banerjee} once again is invoked, and it reduces the problem to showing that for any collection of edges $e_1, \dots, e_q$ of $G$ (not necessarily distinct), we have
$$\reg (I^{q+1} : e_1 \dots e_q) \le \nu+1.$$
To this end, let $J$ be the polarization of $I^{q+1}:e_1 \dots e_q$. It can be seen that
$$J = I(H) + (x_{i_1}y_{i_1}, \dots, x_{i_t}y_{i_t}),$$
where $H$ is a graph over the vertices $\{x_1, \dots, x_n\}$ and $x_{i_1}^2, \dots, x_{i_t}^2$ are the non-squarefree generators of $I^{q+1}: e_1 \dots e_q$. By using standard short exact sequences, it can be shown that
$$\reg(J) = \reg I(H).$$

Observe further that $H$ contains $C_n$ as a subgraph, so $H$ has a Hamiltonian cycle. Thus, the required upper bound now follows from a more general upper bound for the regularity of a graph admitting either a Hamiltonian path or a Hamiltonian cycle (which is the content of \cite[Theorems 3.1 and 3.2]{BHT}).
\end{sketch}

The following theorem was proved in the special case of \emph{whiskered} cycles by Moghimian, Seyed Fakhari, and Yassemi \cite{MFY}, and then in more generality for unicyclic graphs by Alilooee, Beyarslan and Selvaraja \cite{ABS}.

\begin{theorem}[\protect{\cite[Theorem 1.2]{ABS}, \cite[Proposition 1.1]{MFY}}] \label{thm.unicyclic}
Let $G$ be a unicyclic graph and let $I = I(G)$ be its edge ideal. Then for all $q \ge 1$, we have
$$\reg I^q = 2q + \reg I - 2.$$
\end{theorem}

\begin{sketch} The proof is based on establishing the conjectural bound of Conjecture \ref{conj.BBH}.(2)
\begin{align}
\reg (I^q ) \le 2q+\reg I-2. \label{eq.506}
\end{align}

The inequality (\ref{eq.506}) is proved by using induction on $q.$ It also requires a good understanding of $\reg(I(C_n)^q,f_1,\ldots, f_k)$ where $C_n$ is the cycle in $G$ and $f_1, \ldots, f_k$ are the edges of $G$ that are not in $C_n.$ By making use of the general lower bound of Theorem \ref{thm.generallower} and the first main result of \cite{ABS}, namely,
$$\reg I = \left\{ \begin{array}{rl}
\nu(G) + 2 & \text{if } n \equiv 2 (\text{mod } 3) \text{ and } \nu(G \setminus \Gamma(G)) = \nu(G) \\
\nu(G)+ 1 & \text{otherwise,} \end{array} \right.$$
equality is achieved for the latter case. (Here, $\Gamma(G)$ is a well described subset of the vertices in $G$.)  The proof is completed by showing $\reg ( I(G \setminus \Gamma(G))^q) =2q+\nu(G)$ and using inequality  (\ref{eq.500}).
\end{sketch}

\begin{example}
Let $G$ be the graph depicted in Figure \ref{fig:bicyclic}. Macaulay 2 computations show that
$$\reg I(G)=5, \reg I(G)^2=6, \reg I(G)^3=8, \reg I(G)^4=10, \reg I(G)^5=12.$$
Thus, the formula given in Theorem \ref{thm.unicyclic} does not necessarily hold for a graph containing more than 1 cycle.

\begin{figure}[h]
  \includegraphics[width=0.45\linewidth]{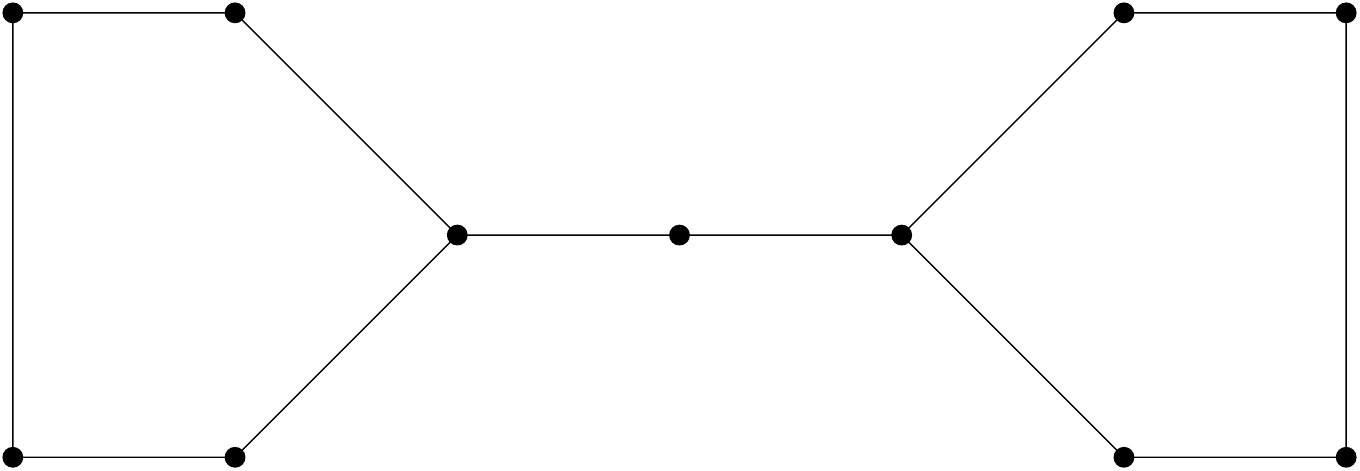}
  \caption{A bicyclic graph}
  \label{fig:bicyclic}
\end{figure}
\end{example}

A particular interesting class of graph is those for which $\nu(G) = 1$. It is expected that for such a graph $G$, powers of its edge ideal should asymptotically have linear resolutions. This is examined in the next two theorems under various additional conditions.

\begin{theorem}[\protect{\cite[Theorem 1.2]{Ba}}] \label{thm.gapcricket}
Let $G$ be a gap-free and cricket-free graph, and let $I = I(G)$ be its edge ideal. Then for all $q \ge 2$, we have
$$\reg I^q = 2q;$$
i.e., $I^q$ has a linear resolution.
\end{theorem}

\begin{sketch} Since $I^q$ is generated in degree $2q$, it remains to show that for $q \ge 2$, $\reg I^q \le 2q$. By induction, making use of the inductive techniques of Theorem \ref{thm.Banerjee} (and \cite[Theorem 3.4]{Ba} which proves that $\reg I \le 3$), it suffices to show that for any collection of edges $e_1, \dots, e_q$ in $G$, $\reg (I^{q+1}:e_1 \dots e_q) \le 2$. Let $J$ be the polarization of $I^{q+1} : e_1 \dots e_q$ and let $H$ be the simple graph associated to $J$. The statement is reduced to showing that $\reg J = 2$, or equivalently (by Theorem \ref{thm.regtwo}), that $H^c$ is a chordal graph. That is, $H$ does not have any induced anticycle of length at least 4.

By contradiction, suppose that $H$ has an anticycle $w_1, \dots, w_s$ of length at least 4. Since an induced anticycle of length 4 gives a gap, we may assume further that $s \ge 5$. Suppose also that $e_1 = xy$. The proof follows from the following observations:
\begin{enumerate}
\item There must be an edge between $\{x,y\}$ and $\{w_1,w_3\}$; otherwise $xy$ and $w_1w_3$ form a gap in $G$.
\item Suppose that $xw_1 \in E_G$. Then neither of $w_2$ nor $w_n$ can coincide with $x$.
\item Neither of $w_2$ nor $w_n$ can coincide with $y$; otherwise the other vertex (among $\{w_2, w_n\}$) and $w_1$ would be even-connected implying that $w_1w_2$ or $w_1w_n$ is an edge in $H$.
\item Either $w_2$ or $w_n$ must be connected to $x$; otherwise, $xy$ and $w_2w_n$ form a gap in $G$.
\item Suppose that $xw_2 \in E_G$. Now, by the same line of arguments applying to $\{w_3, w_n\}$, we deduce that either $w_3$ or $w_n$ must be a neighbor of $x$.
\item Suppose that $xw_3 \in E_G$. We arrive at a contradiction that $\{w_1, w_3, x, y, w_2\}$ forms a cricket in $G$.
\end{enumerate}
\end{sketch}

\begin{theorem}[\protect{\cite[Theorem 4.9]{Erey}}] \label{thm.gapdiamond}
Let $G$ be a gap-free and diamond-free graph, and let $I = I(G)$ be its edge ideal. Then for all $q \ge 2$, we have
$$\reg I^q = 2q;$$
i.e., $I^q$ has a linear resolution.
\end{theorem}

\begin{sketch} The proof of this theorem is based on a good understanding of the combinatorial structures of gap-free and diamond-free graphs. Let $\omega(G)$ denote the largest size of a complete subgraph in $G$ (the \emph{clique number} of $G$). If $\omega(G) < 3$ then $G$ is cricket-free, and the assertion follows from Theorem \ref{thm.gapcricket}. If $\omega(G) \ge 4$, then it is shown that the complement of $G$ is chordal, and the conclusion follows from Theorem \ref{thm.regtwo}.

Consider the case where $\omega(G) = 3$. It is shown that if, in addition, $G$ is also $C_5$-free then the complement of $G$ is either chordal or $C_6$, and the assertion again follows from known results. The essential part of the proof is then to examine the structures of $G$ when $G$ is gap-free, diamond-free and contains a $C_5$. This is where novel and interesting combinatorics happen. It is shown that, in this case, $G$ can be obtained from a list of ten specific graphs via a so-called process of \emph{multiplying vertices}. The proof is completed with a careful analysis of each of these ten graphs to show that $\reg (I(G)^{q+1} : e_1 \dots e_q) \le 2$ and to employ Banerjee's result, Theorem \ref{thm.Banerjee}.
\end{sketch}

Another large class of graphs for which the regularity of powers of their edge ideals can be computed explicitly is that of very well-covered graphs. The following theorem was proved by Norouzi, Seyed Fakhari, and Yassemi \cite{NFY} for very well-covered graphs with an additional condition, and then by Jayanthan and Selvaraja \cite{JS} in full generality for any very well-covered graph.

\begin{theorem}[\protect{\cite[Theorem 5.3]{JS}, \cite[Theorem 3.6]{NFY}}] \label{thm.verywell}
Let $G$ be a very well-covered graph and let $I = I(G)$ be its edge ideal. Then for all $q \ge 1$, we have
$$\reg I^q = 2q + \nu(G)-1.$$
\end{theorem}

\begin{sketch} The proof also uses inductive techniques given by Theorem \ref{thm.Banerjee} and goes along the same line as that of previous theorems in this section. The heart of the arguments is to verify that for any collection of edges $e_1, \dots, e_q$ in $G$, we have
$$\reg (I^{q+1} : e_1 \dots e_q) \le \nu(G)+1.$$

If $J = (I^{q+1}: e_1 \dots e_q)$ is squarefree then this is achieved by letting $H$ be the graph associated to $J$ and establishing the following facts:
\begin{enumerate}
\item $H$ is also very well-covered;
\item $\nu(H) \le \nu(G)$;
\item $\reg I(H) = \nu(H) + 1$ and $\reg I(G) = \nu(G)+1$ (this is the content of \cite[Theorem 4.12]{MMCRTY}).
\end{enumerate}

The arguments are much more involved for the case where $J$ is not squarefree. Let $H$ be the graph corresponding to the polarization of $J$. The proof is completed by an ingenious induction on the number of vertices added to $G$ in order to obtain $H$.
\end{sketch}

\section{Higher Dimension} \label{chapter6}

In this section we discuss the regularity of  squarefree monomial ideals generated in degree more than two. The goal, as in the case of the edge ideals, is to find bounds and/or formulas for regularity of the ideal in terms of ``its combinatorics". As mentioned in the preliminaries, one can view these ideals both as Stanley-Reisner ideals or edge ideals of hypergraphs. If interpreted as Stanley-Reisner ideals, via Hochster's formula one can potentially get all Betti numbers of a given ideal in terms of the combinatorics of the underlying complex. The hypergraph case, however, is far less understood compared to either the edge ideal case or the Stanley-Reisner ideal interpretation. However, in the last decade, some results have been proven in that direction.  In some cases, general squarefree monomial ideals can be viewed as path ideals of simple graphs, and there has been some progress for those cases in the last few years.

 Since the main topic of this survey is the edge ideals, we are not trying to be comprehensive in this section. Our aim is to give the reader an idea of possible approaches to generalize the results on edge ideals to higher dimensions. We split the section into three subsections; devoted to Stanley-Reisner ideals, edge ideals of hypergraphs and path ideals respectively.
\subsection{Stanley-Reisner Ideals}

   As mentioned earlier, any squarefree monomial ideal can be viewed as the Stanley-Reisner ideal of some simplicial complex. This combinatorial interpretation is by far the most studied one among monomial ideals generated in higher degrees. %The Hochster's formula mentioned in section three remains the pivotal result. It is worth remembering that via polarization, any monomial ideal can be converted to a squarefree monomial ideal. So any result regarding the regularity of powers of edge ideal can be viewed as result regarding regularity of some Stanley-Reisner ideal. However computing that simplicial complex from the underlying graph is in general extremely hard.  In this subsection we mention two results, one classification and one upper bound.

    The first result we mention is a famous theorem by Eagon and Reiner which establishes the relation between minimal regularity and maximal depth as well as relation between regularity of a Stanley-Reisner ideal and the structure of a special kind of dual complex, namely, the Alexander dual. In a sense it is similar to Fr\"oberg's result; this gives a classification of linear resolution. The\textit{ Alexander dual} of a simplicial complex $\Delta$ is the simplicial complex whose faces are the complements of the nonfaces of $\Delta$. If $I$ is a Stanley-Reisner ideal then by $I^\vee$ we denote the Stanley-Reisner ideal of the Alexander dual of the simplicial complex of $I$.

\begin{theorem}[\protect{\cite[Theorem 3]{ER}}]
Let $I$ be a Stanley-Reisner ideal in $S=K[x_1,\dots,x_n]$. Then $I$ has $q$-linear resolution if and only if depth of $S/I^\vee$ has depth $n-q$. In particular, $I$ has linear resolution if and only if the Alexander dual of $\Delta (I)$ is Cohen-Macaulay.
\end{theorem}

  The next result bounds regularity with another important algebraic invariant, the arithmetic degree. For a squarefree monomial ideal $I$, its \textit{arithmetic degree}, denoted by $\adeg(I)$, is given by number of facets in the corresponding simplicial complex.

\begin{theorem}[\protect{\cite[Theorem 3.8]{FkT}}]
Suppose $I$ is a Stanley-Reisner ideal with codimension $\geq 2$. Then, $\reg I \le \adeg(I)$.
\end{theorem}

\subsection{Hypergraphs}

   Any squarefree monomial ideal is the edge ideal of a  hypergraph. The corresponding simplicial complex is the independence complex of this hypergraph and computing that is an NP-hard problem in general.%  For splittable monomial ideals some interesting results regarding regularity has been proved.

%  A monomial ideal $I$ is said to be \textit{splittable} if there exists two nonzero monomial ideals $J,K$ such that $I=J+K$ where the set of minimal monomial generators $T$ of $I$ is a disjoint union of those of $J$ and $K$ ($T_1$ and $T_2$ respectively) and if there is a function $f: T \longrightarrow T_1 \times T_2$ satisfying, for all $a\in T$ and $A\subseteq T$, $a$ is the lcm of coordinates of $f(a)$ and lcm of both first and second coordinates of $f(A)$ strictly divide $f(A)$.

%\begin{theorem}[H\`a, Van-Tuyl]
%Let $H$ be a hypergraph with at least two edges. An edge $m$ is a splitting edge if and only if there exists a vertex $z\in E(H)$ such that $$ <m> \cap I(H \setminus V(m)) \subseteq <m> \cap I(H \setminus \{z\})$$ Here $H \setminus \{z\}$ denotes the hypergraph obtained from $H$ by deleting all edges containing $z$.
%\end{theorem}

%\begin{theorem}[H\`a, Van-Tuyl]
%Let $H$ be a properly connected hypergraph and let $m$ be a splitting edge of $H$. Let $d$ be the degree of $m$, $A$ be the set of all edges of $H$ with distance greater than or equal to $d+1$ from $m$ and $t$ be the cardinality of $N(m)$. Then for all $i\geq 1$ we have the following $$\beta_{ij}( I(H)) = \beta_{ij}( I(H \setminus V(m))) + \Sigma \frac{t!}{l!t-l!} \beta_{i-1-l j-d-l} I(A)$$
%\end{theorem}

The problem of finding bounds for the regularity of graphs can be extended to hypergraphs. In \cite{CHHKTT} the authors provided a sufficient condition for a hypergraph to have regularity $\leq 3.$ For every vertex $x$ of a hypergraph $\H$, let $\H : x$ denote the simple hypergraph of all minimal subsets $A \subset  V \setminus \{x\}$  such that  $A$ or $A\cup \{x\}$ is an edge of $\H.$

\begin{theorem}[\protect{\cite[Theorem 6.4]{CHHKTT}}] \label{thm.hyperbound}
Let $\H$ be a simple hypergraph such that $\H :x$ is a graph whose complement is chordal for all vertices of $\H.$ Then, $ \reg \H \leq 3.$
\end{theorem}

Any hypergraph whose all edges have same cardinality is called \textit{uniform}. A $d$-uniform hypergraph (i.e., a uniform hypergraph whose all edges are represented by monomials of degree $d$) is called \textit{properly-connected} if for any two edges $E,E'$ sharing at least one vertex, the length of the shortest path between $E$ and $E'$ is $d- |E \cap E'|$. If the length of the shortest path between two edges is $\geq t$ then they are called $t$-disjoint (for relevant definitions see \cite{HVT2008}).

\begin{example}
Consider the $4$-uniform hypergraph H with edge set: $$\{x_1x_2x_3x_4, x_1x_2x_3x_7, x_1x_2x_6x_7, x_1x_5x_6x_7, x_1x_5x_6x_8\}$$ There is a proper irredundant chain of length $4$ from the edge $E=x_1x_2x_3x_4$ to $E_0=x_1x_5x_6x_8$. Furthermore, there is no shorter such chain. But these edges have a nonempty intersection. So H is not properly-connected since the distance between them is $4$ which is not same as $4 -| E \cap  E_0| = 3$. This hypergraph is not properly-connected. On the other hand every finite simple graph is properly connected.
\end{example}

A $d$-uniform properly connected hypergraph $H$ is called \textit{triangulated} if for any subset $A$ of vertices the induced subgraph $H'$ on $A$ has a vertex $x$ such that the induced subgraph on its neighborhood $N(x) \cup \{x\}$ is complete. For a more formal definition and relevant notions see \cite{HVT2008}.

\begin{example}
Simple graphs that are triangulated are precisely the chordal graphs due to Theorem 5.2 of \cite{HVT2008}.
\end{example}

The following result by H\`a and Van-Tuyl gives a lower bound for regularity of properly connected hypergraphs and the bound becomes an equality if the hypergraph is triangulated.

\begin{theorem}[\protect{\cite[Theorem 6.5 and 6.8]{HVT2008}}]
Let $H$ be a properly-connected hypergraph. Suppose $d$ is the common cardinality of the edges in $H$. Let $c$ be the maximal number of pairwise $d+1$ disjoint edges in $H$. Then, $\reg H \ge (d-1)c + 1$, and the equality occurs if $H$ is triangulated.
\end{theorem}

A  $2$-\textit{collage} for a hypergraph is a subset $C$ of the edges with the property that for every edge $E$  of the hypergraph, we can delete a vertex $v$ of $E$ so that $E \setminus \{v\}$ is contained in some edge of $C$.  For uniform hypergraphs, the condition for a collection $C$ of the edges to be a $2$-collage is equivalent to requiring that for any edge $E$ not in $C$, there exists  $F\in C$ such that the symmetric difference of  $E$ and $F$ consists of exactly two vertices. When H is a graph, it is straightforward to see that for any minimal $2$-collage, there is a maximal matching of the same or lesser cardinality.

The following theorem gives a formula for regularity of uniform hypergraphs in terms of its collages.

\begin{theorem}[\protect{\cite[ Theorem 1.1]{HW}}]
Let $H$ be a simple $d$-uniform hypergraph and let $c$ be the minimum size of a $2$ collage in $H$. Then, $\reg I(H) \leq (d-1)c+1$.
\end{theorem}

The next theorem is a generalization of the previous that work for all simple hypergraph.
\begin{theorem}[\protect{\cite[Theorem 1.2] {HW}}]
Let $H$ be a simple hypergraph and let $\{m_1,....,m_c\}$ be a $2$-collage. Then, $\reg I(H) \leq \sum_{i=1}^c |m_i| - c + 1$.
\end{theorem}

The next theorem gives a formula for regularity of edge ideals of \textit{clutters}, which are hypergraphs where no edge contains any other edge as a subset. A monomial ideal $I$ has \textit{linear quotients } if the monomials that generate $I$ can be ordered $g_1, \ldots, g_q$ such that for all $1\leq i \leq q-1, ((g_1, \ldots, g_i):g_{i+1})$ is generated by linear forms. For further details see \cite{MV}.

\begin{theorem}[\protect{\cite[Corollary 3.35]{MV}}]
 Let $H$ be a clutter such that $I(H)$ has linear quotients. Then, $\reg I(H) =\max\{|E|: E \in E(H)\}$.
\end{theorem}

Taylor resolutions, introduced in \cite{Tay}, are free resolutions of monomial ideals constructed in a natural way that makes computations easy to deal with. These are in general not minimal. The following theorem gives a formula for regularity when Taylor resolutions are minimal.  A hypergraph $H$ is called \textit{saturated} when the Taylor resolution of $I(H)$ is minimal.  For further details see \cite{LM}.

\begin{theorem}[\protect{\cite[Proposition 4.1]{LM}}]
For a saturated hypergraph $H$, we have $\reg I(H =|X|-|L| +1$. Here $L$ is the minimal monomial generating set of $I(H)$ and $X$ is the number of variables that divide some minimal monomial generator.
\end{theorem}

\subsection{Path Ideal}

    Path ideals of finite simple graphs are interesting generalizations of edge ideals, see Section \ref{sub.alg} to recall the definition of path ideals. For gap-free graphs their regularity  too tend to behave nicely.

\begin{theorem}[\protect{\cite[Theorem 1.1]{Bp}}]
Let $G$ be a gap-free and claw-free graph. Then for all $t \leq 6$, $I_t$ has linear resolution. If further $G$ is induced whiskered-$K_4$ free then for all $t$, $I_t$ has linear resolution.
\end{theorem}

Let $L_n$ denote a line of length $n$. The regularity of path ideals of lines have been computed by Alilooee and Faridi in the next theorem.

\begin{theorem}[\protect{\cite[Theorem 3.2]{AlFr}}]
Let $n,t,p$ and $d$ be integers such that $n \geq 2$, $2 \leq t \leq n, n=(t+1)p+d$, where $p\geq 0, 0 \leq d \leq t$. Then, $\reg I_t (L_n)$ is $p(t-1)+1$ for $d < t$ and $(p+1)(t-1) +1$ for $d=t$.
\end{theorem}

We end this section by a result regarding a somewhat different kind of path ideals. A tree is a graph in which there exists a unique path between every pair of distinct vertices; a rooted tree is a tree together with a fixed vertex called the root. In particular, in a rooted tree there exists a unique path from the root to any given vertices. We can also view a rooted tree as a directed graph by assigning to each edge the direction that goes ``away'' from the root.  Let $\Gamma$ be a rooted tree and $I_t(\Gamma)$ be the squarefree monomial ideal generated by all ``directed paths'' of length $t-1$ in the above sense. The following theorem gives regularity of such ideals.

In this theorem, we define $l_t(\Gamma)$ to be the number of leaves in $\Gamma$ whose level is at least $t -1$ and $p_t(\Gamma)$ to be the maximal number of pairwise disjoint paths of length $t$ in $\Gamma$ (i.e., $p_t(\Gamma)= \max \{|D|~|~ D \text{ is a set of disjoint paths of length } t \text{ in } \Gamma\}$).

\begin{theorem}[\protect{\cite[Theorem 3.4]{BHO}}]
Let $\Gamma$ be a rooted tree on $n$ vertices. Then, $\reg I_t(\Gamma) \leq (t-1)[l_t(\Gamma) + p_t(\Gamma)] +1$.
\end{theorem}

\section{Open problems and questions} \label{chapter7}

We end the paper by listing a number of open problems and questions in the research area which we hope to be solved.

Our first problem is inspired by Theorems \ref{thm.regtwo} and \ref{thm.regbipartitethree}. Since graphs with regularity 2 are classified. The next class of graph to examine is that of regularity 3.

\begin{problem} \label{prob.reg3} 
Characterize graphs $G$ for which $\reg I(G) = 3$.
\end{problem}

In various results, for example Theorems \ref{thm.gapclawfree} and \ref{thm.hyperbound}, ``local'' condition on the regularity of $G:x$ for all vertices $x$ lead to a ``global'' statement on the regularity of $G$. We ask if similar local conditions on $G:x$, for all vertices $x$, would also lead to a statement on the asymptotic linear function $\reg I(G)^q$.

\begin{question} Let $G$ be a simple graph and let $I = I(G)$ be its edge ideal. Suppose that for any vertex $x$ in $G$, we have $\reg (I:x) \le r$. Does this imply that for any $q \ge 1$,
$$\reg I^q \le 2q + r-1?$$
\end{question}

As noted in Lemma \ref{lem.deletecontract}, for any vertex $x$, the regularity of $I(G)$ is always equal to either the regularity of $I(G) : x$ or the regularity of $(I(G),x)$. It would be interesting to know for which vertex $x$ the equality happens one way or another.

\begin{problem} Let $G$ be a graph and let $I = I(G)$ be its edge ideal. Find conditions on a vertex $x$ of $G$ such that
\begin{enumerate}
\item $\reg I = \reg (I,x).$
\item $\reg I = \reg (I:x).$
\end{enumerate}
\end{problem}

The regularity of $I(G)$ has been computed for several special classes of graphs (see Theorem \ref{thm.inducedmatch1}). A particular class of graphs which is of interest is that of vertex decomposable graphs. For a vertex decomposable graph $G$, the statement of Lemma \ref{lem.deletecontract} can be made slightly more precise (see \cite{HW}), namely, there exists a vertex $x$ such that
$$\reg G = \max \{\reg (G \setminus N_G[x]) +1, \reg (G \setminus x)\}.$$
For such a graph $G$, it is also known that the independent complex $\Delta(G)$ is shellable and the quotient ring $S/I(G)$ is sequentially Cohen-Macaulay.

\begin{problem} Let $G$ be a vertex decomposable graph. Compute $\reg I(G)$ via combinatorial invariant of $G$.
\end{problem}

As noted throughout Sections \ref{chapter4} and \ref{chapter5}, the induced matching number of a graph $G$ is closely related to the regularity of $I(G)$. In fact, $\nu(G) + 1$ gives a lower bound for $\reg I(G)$ and, more generally, $2q+ \nu(G)-1$ gives a lower bound for $\reg I(G)^q$ for any $q \ge 1$. Moreover, for many special classes of graphs the equality has been shown to hold. Thus, it is desirable to characterize all graphs for which the equality is attained.

\begin{problem} \label{prob.nu}
Characterize graphs $G$ for which the edge ideals $I = I(G)$ satisfy
\begin{enumerate}
\item $\reg I = \nu(G)+1$.
\item $\reg I^q = 2q+\nu(G) -1 \ \forall \ q \gg 0.$
%\item $\reg I^q = 2q+\reg(I) - 2 \ \forall \ q \gg 0.$
\end{enumerate}
\end{problem}

When $\nu(G) = 1$, the answer to Problem \ref{prob.nu}.(2) is predicted in the following open problem.

\begin{problem}[Francisco-H\`a-Van Tuyl and Nevo-Peeva] \label{prob.FHVT}
Suppose that $\nu(G) = 1$, i.e., $G^c$ has no induced 4-cycle and let $I = I(G)$.
\begin{enumerate}
\item Prove (or disprove) that $\reg I^q = 2q$ for all $q \gg 0$.
\item Prove (or disprove) that $\reg I^{q+1} = \reg I^q + 2$ for all $q \ge \reg(I)-1$.
\end{enumerate}
\end{problem}

Note that examples exist in which $\nu(G) = 1$ and $\reg I^q \not= 2q$ for small values of $q$ (cf. \cite{NP}), so in Problem \ref{prob.FHVT} it is necessary to consider $q \gg 0$. A satisfactory solution to Problem \ref{prob.reg3} would be a good starting point to tackle Problem \ref{prob.FHVT}, since regularity 3 is the first open case of the problem. In fact, in this case, $\reg I(G)^q$ is expected to be linear starting at $q = 2$.

\begin{problem}[Nevo-Peeva]
Suppose that $\nu(G) = 1$ and $\reg I(G) = 3$. Then is it true that for all $q \ge 2$,
$$\reg I(G)^q = 2q?$$
\end{problem}

In many known cases where the asymptotic linear function $\reg I(G)^q$ can be computed, it happens to be $\reg I(G)^q = 2q+\nu(G)-1$. We would like to see if this is the case when the equality is already known to hold for small values of $q$. If this is indeed the case then how far one must go before concluding that the equality holds for all $q \ge 1$?

\begin{problem} Let $G$ be a graph. Find a number $N$ such that if $\reg I(G)^q = 2q+ \nu(G)-1$ for all $1 \le q \le N$ then, for all $q \ge 1$, we have
$$\reg I(G)^q = 2q+\nu(G)-1.$$
\end{problem}

\noindent (Computational experiments seem to suggest that $N$ can be taken to be 2.)

A particularly related question is whether for special classes of graphs satisfying $\reg I(G) = \nu(G)+1$ one would have $\reg I(G)^q = 2q+\nu(G)-1$ for all $q \ge 1$ (or for all $q \gg 0$). Inspired by Theorem \ref{thm.inducedmatch1}, we raise the following question.

\begin{question} \label{quest.nu+1}
Let $G$ be a graph and let $I = I(G)$ be its edge ideal. Suppose that $G$ is of one of the following types:
\begin{enumerate}
\item $G$ is chordal;
\item $G$ is weakly chordal;
\item $G$ is sequentially Cohen-Macaulay bipartite;
\item $G$ is vertex decomposable and contains no closed circuit of length 5;
\item $G$ is $(C_4,C_5)$-free vertex decomposable.
\end{enumerate}
Is it true that for all $q \gg 0$,
$$\reg I^q = 2q+\nu(G)-1?$$
\end{question}

Note that chordal graphs and sequentially Cohen-Macaulay bipartite graphs are vertex decomposable. Also bipartite graphs contain no closed circuit of length 5. Thus, in Question \ref{quest.nu+1}, an affirmative answer to part (4) would imply that for part (3).

It is well-known that the resolution of a monomial ideal is dependent on the characteristic of the ground field. And yet, in all known cases where the regularity of powers of an edge ideal can be computed, it is characteristic-independent. We would like to see examples where this is not the case, or a confirmation that this is always the case if the regularity of the edge ideal itself is characteristic-independent.

\begin{problem} \quad
\begin{enumerate}
\item Find examples of graphs $G$ for which the asymptotic linear function $\reg I(G)^q$, for $q \gg 0$, is characteristic-dependent.
\item Suppose that $\reg I(G)$ is independent of the characteristic of the ground field. Is the asymptotic linear function $\reg I(G)^q$, for $q \gg 0$, necessarily also characteristic-independent?
\end{enumerate}
\end{problem}

When an exact formula may not be available, it is of interest to find a lower and an upper bound. Since the lower bound in Theorem \ref{thm.generallower} holds for any graph, it is desirable to find an upper bound to couple with this general lower bound. In this direction, one may either try to prove the bound in Theorem \ref{thm.JNS} for any graph or to relate the regularity of powers of $I(G)$ to the regularity of $I(G)$ itself.

\begin{conjecture}[Alilooee, Banerjee, Beyarslan and H\`a] \label{conj.BBH}
Let $G$ be any graph and let $I = I(G)$ be its edge ideal. Then for all $q \ge 1$, we have
\begin{enumerate}
\item $\reg I^q \le 2q + \text{co-chord}(G) -1.$
\item $\reg I^q \le 2q+\reg I -2.$
\item $\reg I^{q+1} \le \reg I^q + 2.$
\end{enumerate}
\end{conjecture}

Note that by Theorem \ref{thm.reggeneralupper}, Conjecture \ref{conj.BBH}.(3) $\Rightarrow$ Conjecture \ref{conj.BBH}.(2) $\Rightarrow$ Conjecture \ref{conj.BBH}.(1).

Even though, in general, $\reg I^q$ is asymptotically a linear function, for small values of $q$, there are examples for which $\reg I^q > \reg I^{q+1}$. We ask if this would \emph{not} be the case for edge ideals of graphs.

\begin{question} Let $G$ be a graph and let $I = I(G)$ be its edge ideal. Is the function $\reg I^q$ increasing for all $q \ge 1$?
\end{question}

In investigating the asymptotic linear function $\reg I^q$, it is also of interest to know the smallest value $q_0$ starting from which $\reg I^q$ attains its linear form. In known cases where the regularity of $I(G)^q$ can be computed explicitly for all $q \ge 1$, we have $q_0 \le 2$. Computational experiments seem to suggest that this is indeed always the case.

\begin{question} \label{quest.q0}
Let $G$ be a graph and let $q_0$ be the least integer such that $\reg I(G)^q$ is a linear function for all $q \ge q_0$. Is it true that $q_0 \le 2$?
\end{question}

A much weaker question, yet still very interesting, would be whether $q_0 \le \reg I(G)$? --- This question is still wide open! A significant step toward answering Question \ref{quest.q0} is to give a uniform bound for $q_0$ (i.e., a bound that does not depend on $G$).

\begin{problem} Find a number $N$ (which may depend on $n$ and $m$) such that for any graph $G$ over $n$ vertices and $m$ edges, we have $q_0 \le N$.
\end{problem}

\end{document}